\documentclass[11pt]{article}
\usepackage{latexsym}
\usepackage{amsfonts,amssymb,amsmath}

\usepackage[yyyymmdd]{datetime}

\newtheorem{thm}{Theorem}[section]
\newtheorem{cor}[thm]{Corollary}
\newtheorem{lem}[thm]{Lemma}
\newtheorem{pro}[thm]{Proposition}
\newtheorem{defn}[thm]{Definition}

\newcommand{\ov }{\overline }
\newcommand{\mal}{{\footnotesize \textcircled{{\scriptsize \sf m}}\ }}
\newcommand{\liteq}{\stackrel{+}{=}}

\title{Global local covers}

\author{John L.\ Rhodes, Benjamin Steinberg, J.C.\ Birget}
 
\date{\scriptsize{
2.iv.2019}}

\begin{document}
\maketitle

\begin{abstract}
This paper gives a systematic construction of certain covers of finite 
semigroups. These covers will be used in future work on the complexity 
of finite semigroups.
\end{abstract}



\section{Introduction}

This paper is an incomplete description of the important {\em glc cover} 
for finite semigroups, which will play a vital role in the future proposed 
proof of the decidability of {\em c}. 
(See see e.g.\ \cite{EilenbergB}, \cite{qTheory} for the definition of the 
complexity {\em c} of finite semigroups.)

\smallskip 

In the following, certain surmorphisms between semigroups, or between
automata, will be called {\em covers}. 
Covers are surmorphisms, but they are obtained by constructing a 
pre-image in a ``systematic way''. Expansions (as in \cite{Afinexp, 
ItExp}) are examples of covers; in general, covers are a weaker form of 
expansions.


\section{A direct construction of the glc-cover}

\subsection{Lattices of automata}

We consider the categories ${\mathcal S}_A$ and ${\mathcal M}_A$ of finite
semigroups, respectively monoids, over a chosen generating set $A$; the
morphisms are the usual semigroup (or monoid) morphisms that commute with 
the map that sends $A$ to a generating set of the semigroup (or monoid). 
A semigroup in that category will usually be denoted by $S_A$ or $(S,A)$.  

We consider also the category ${\mathcal A}_A$ of finite automata with 
alphabet $A$. Here, the objects are all deterministic finite automata with 
input alphabet $A$, with a chosen start state {\sc i}, that are trim (i.e., 
all states are reachable from {\sc i}), and such that {\sc i} cannot be
returned to. The morphisms in ${\mathcal A}_A$ are all state maps that commute 
with the next-state function and preserve the start state.
Usually, we work in the category of {\it isomorphism classes} of
${\mathcal A}_A$, which we also denote by ${\mathcal A}_A$; but because of
the role of the start state this makes no difference. 

For any semigroup $S$ generated by $A$, the {\it right Cayley graph} of 
$S$ belongs to ${\mathcal A}_A$. For clarity, let us  define the right
Cayley graph $\Gamma_A(S, \sigma)$ for a semigroup $S$ and a map 
$\sigma: A \to S$ that maps $A$ onto a generating set of 
$S$: $\Gamma_A(S, \sigma)$ is a directed labeled rooted graph, with 
vertex set $S^I$ (i.e., $S$ with a new identity 
element $I$ added, even in case $S$ already had an identity), with root $I$, 
and with edge set $\{ (s, a, s \cdot \sigma(a)) : s \in S^I, \, a \in A\}$ 
(where $\cdot$ is the product in $S$). The directed edge 
$(s, a, s \cdot \sigma(a))$ is labeled by $a$. We often write $\Gamma(S_A)$
instead of $\Gamma_A(S, \sigma)$. 

For automata ${\bf A}_1, {\bf A}_2$ in ${\mathcal A}_A$ we write 
${\bf A}_1 \twoheadleftarrow {\bf A}_2$ iff there exists a surjective 
${\mathcal A}_A$-morphism from ${\bf A}_2$ onto ${\bf A}_1$.
The relation $\twoheadleftarrow$ is a pre-order on ${\mathcal A}_A$.
Because of the role of the start state, a surjective morphism is unique if 
it exists (for a given pair ${\bf A}_1, {\bf A}_2$); hence, this pre-order 
is a partial order.  We also write ${\bf A}_2 \twoheadrightarrow {\bf A}_1$ 
for ${\bf A}_1 \twoheadleftarrow {\bf A}_2$.

For automata ${\bf A}_1, {\bf A}_2$ in ${\mathcal A}_A$ we define the
{\it direct product} in the category ${\mathcal A}_A$. 
For $i = 1,2$, let ${\bf A}_i = (Q_i, I_i, \cdot_i)$, where 
$\cdot_i : Q_i \times A \to Q_i$ is the next-state function. Then

\smallskip

${\bf A}_1 \times_A {\bf A}_2$ $ \ = \ $
$\big( \{(I_1 \cdot_1 w, \ I_2 \cdot_2 w) \in Q_1 \times Q_2 : $
$w \in A^*\}, \ (I_1,I_2), \ \cdot \big)$, 

\smallskip

\noindent where for all $w \in A^*$, $a \in A$, we define 
 \ $(I_1 \cdot_1 w, \ I_2 \cdot_2 w) \cdot a$ $\ = \ $ 
$(I_1 \cdot_1 w \cdot_1 a, \ I_2 \cdot_2 w \cdot_2 a)$.  

For two automaton morphisms $\varphi_i: {\bf A}_i \to {\bf B}$ in the
category ${\mathcal A}_A$ we can define the {\it direct product}
$\varphi_1 \times_A \varphi_2: {\bf A}_1 \times_A {\bf A}_2 \to {\bf B}$
in ${\mathcal A}_A$ by
 \ $(I_1 \cdot_1 w, I_2 \cdot_2 w)\varphi_1 \times_A \varphi_2$
$\ = \ I \cdot_{_{\bf B}} w$ \ (for all $w \in A^*$).
This is well-defined: If \ $(I_1 \cdot_1 u, \ I_2 \cdot_2 u)$
$=$ $(I_1 \cdot_1 v, \ I_2 \cdot_2 v)$ \ then
 \ $I_i \cdot_i u = I_i \cdot_i v$ \ (for $i = 1, 2$), hence
$(I_i \cdot_i u)\varphi_i = (I_i \cdot_i v)\varphi_i$. Moreover,
$(I_i \cdot_i u)\varphi_i = (I_i)\varphi_i \cdot_{_{\bf B}} u$
$ = I \cdot_{_{\bf B}} u$; similarly,
$(I_i \cdot_i v)\varphi_i = I \cdot_{_{\bf B}} v$. Thus,
$I \cdot_{_{\bf B}} u = I \cdot_{_{\bf B}} v$.

We show next that $\twoheadleftarrow$ is a {\it lattice} order. 

\begin{lem}
(1)  The {\em join} ${\bf A}_1 \vee {\bf A}_2$ in ${\mathcal A}_A$ is 
isomorphic to the direct product ${\bf A}_1 \times_A {\bf A}_2$. 
In the case of right Cayley graphs of semigroups in ${\mathcal S}_A$, the 
join is (up to isomorphism) the direct product in ${\mathcal S}_A$.

\noindent (2)  The {\em meet} is determined by the join as  
 \ ${\bf A_1} \wedge {\bf A}_2 \ = \ \bigvee $
$\{ {\bf A} : {\bf A}_1 \twoheadrightarrow {\bf A}$
and ${\bf A}_2 \twoheadrightarrow {\bf A} \}$ \ (up to isomorphism).  
\end{lem}
{\bf Proof.} 
(1) The projections from ${\bf A}_1 \times_A {\bf A}_2$ onto ${\bf A_1}$,
respectively ${\bf A}_2$, show that ${\bf A}_1 \times_A {\bf A}_2$ is an 
upper bound on ${\bf A_1}$ and ${\bf A}_2$ for the order $\twoheadleftarrow$. 
To show that it is the least upper bound, we consider the commutative 
diagrams in the category ${\mathcal A}_A$ that defines the direct product, 
with arrows 
${\bf A}_1 \twoheadleftarrow  {\bf X} \twoheadrightarrow {\bf A}_2$, 
${\bf X} \twoheadrightarrow {\bf A}_1 \times_A {\bf A}_2$, and the 
projections ${\bf A_1} \twoheadleftarrow {\bf A}_1 \times_A {\bf A}_2$ 
$\twoheadrightarrow {\bf A}_2$, for any object ${\bf X}$ in ${\mathcal A}_A$.
The join is defined by exactly the same diagrams, so the join is the direct 
product.

For Cayley graphs $\Gamma_A(S_1, \sigma_1)$ and $\Gamma_A(S_2, \sigma_2)$, 
viewed as $A$-automata, the direct product 
$\Gamma_A(S_1, \sigma_1) \times_A \Gamma_A(S_2, \sigma_2)$ has the set 
$\{(\sigma_1(w), \sigma_2(w)) : w \in A^* \}$ as its states (vertices); 
here we extended $\sigma_1$ and $\sigma_2$ to homomorphisms from $A^*$ onto
$S_1$ or $S_2$. Hence $\Gamma_A(S_1, \sigma_1) \times_A \Gamma_A(S_2,
\sigma_2)$, with start state $(I_1, I_2)$, is the Cayley graph of 
$(S_1)_A \times_A (S_2)_A$ (semigroup direct product in the category 
${\mathcal S}_A$).

\smallskip

\noindent (2) This follows from the general fact that if a partial order 
has finite joins and a global minimum, and if every principal ideal is 
finite, then this partial order is a lattice. This finiteness comes from 
the fact that here all automata are finite. 
 \ \ \ $\Box$

\medskip  

In the partial order $\twoheadleftarrow$ we define an interval 
 \ $[{\bf A}_1, {\bf A}_2]$ $\ = \ $ $\{ {\bf A} \in {\mathcal A}_A :$
${\bf A}_1 \twoheadleftarrow {\bf A} \twoheadleftarrow {\bf A}_2 \}$.
This interval is a {\it finite lattice} with a global 
maximum and a global minimum. Finiteness comes from the fact that here all
automata are finite.
 
For an automaton ${\bf A}$ in ${\mathcal A}_A$ and for two states 
$q_1, q_2$ of ${\bf A}$, we say that $q_2$ is {\it reachable} from $q_1$
(and we write $q_2 \leq_{\mathcal R} q_1$) iff $q_2 = q_1 \cdot w$
for some $w \in A^*$. This relation is a pre-order, and we write 
$q_2 \equiv_{\mathcal R} q_1$ iff $(q_2 \leq_{\mathcal R} q_1$ and
$q_1 \leq_{\mathcal R} q_2)$.  An $\equiv_{\mathcal R}$ equivalence 
class is called a {\it reachability class}, or an $\mathcal R$-class; 
this is the set of vertices of a strongly connected component of the 
underlying digraph of the automaton.  For two states $q_1, q_2$ we will 
also use the notation $q_2 \neq \equiv_{\mathcal R} q_1$ as a shorthand
for $(q_2 \equiv_{\mathcal R} q_1$ and $q_2 \neq q_1)$.

An automaton morphism $\varphi: {\bf A}_1 \to {\bf A}_2$ is said to be
1:1${\mathcal R}$ iff the restriction of $\varphi$ to any reachability
class of ${\bf A}_1$ is injective. 

\begin{lem}
If two morphisms $\varphi_1: {\bf A}_1 \to {\bf B}$ and
$\varphi_2: {\bf A}_2 \to {\bf B}$ are 1:1${\mathcal R}$ then the direct
product morphism
$\varphi_1 \times_A \varphi_2: {\bf A}_1 \times_A {\bf A}_2 \to {\bf B}$
in ${\mathcal A}_A$ is also 1:1${\mathcal R}$.
\end{lem}
{\bf Proof.}  The state set of ${\bf A}_1 \times_A {\bf A}_2$ is
$\{(I_1 \cdot_1 w, I_2 \cdot_2 w) : w \in A^*\}$.  
The morphism $\varphi_1 \times_A \varphi_2$ is defined by
 \ $\varphi_1 \times_A \varphi_2(I_1 \cdot_1 w, I_2 \cdot_2 w)$
$ = I \cdot_{_{\bf B}} w$
$= \varphi_1(I_1 \cdot_1 w) = \varphi_2(I_2 \cdot_2 w)$. The last two 
equalities hold because $\varphi_1$ and $\varphi_2$ are morphisms; in 
particular, $\varphi_1(I_1) = I = \varphi_2(I_2)$.  If \   
$(I_1 \cdot_1 u, I_2 \cdot_2 u)$ $\neq \equiv_{\mathcal R}$
$(I_1 \cdot_1 v, I_2 \cdot_2 v)$ \ then
 \ $I_1 \cdot_1 u \equiv_{\mathcal R} I_1 \cdot_1 v$ \ and \  
$I_2 \cdot_2 u \equiv_{\mathcal R} I_2 \cdot_2 v$; moreover, 
 \ $I_1 \cdot_1 u  \neq I_1 \cdot_1 v$ \ or \ 
 $I_2 \cdot_2 u  \neq I_2 \cdot_2 v$.
If \ $I_1 \cdot_1 u  \neq I_1 \cdot_2 v$ \ then the 1:1${\mathcal R}$ 
property of $\varphi_1$ implies 
 \ $\varphi_1(I_1 \cdot_1 u) \neq \varphi_1(I_1 \cdot_1 v)$,
hence \ $\varphi_1 \times_A \varphi_2(I_1 \cdot_1 u, I_2 \cdot_2 u)$ 
$\neq$ $\varphi_1 \times_A \varphi_2(I_1 \cdot_1 v, I_2 \cdot_2 v)$.
If \ $I_2 \cdot_2 u  \neq I_2 \cdot_2 v$ \ the same conclusion holds.
Thus, $\varphi_1 \times_A \varphi_2$ is 1:1${\mathcal R}$.
 \ \ \ $\Box$

\medskip

For a finite semigroup $S$ and $x \in S$, let $x^{\omega}$ be the idempotent 
in $\{x^n : n > 0\}$. The {\it algebraic rank} of $x$ is defined to be the 
length of a longest strict $\mathcal J$-chain of {\em regular}
$\mathcal J$-classes, ascending from $x^{\omega}$ (see \cite{GST}). 
The rank of $x$ is a non-negative integer (possibly 0); we denote it by 
${\sf rank}_S(x)$ or ${\sf rank}(x)$. Formally,  

\medskip

${\sf rank}_S(x) \ = \ {\sf max}\{r \ : \ x^{\omega} <_{\mathcal J} J_1$
    $<_{\mathcal J} \ \ldots \ <_{\mathcal J}  J_r,$ \ where 
    $J_1, \ldots, J_r$ are regular $\mathcal J$-classes of $S \}$.
 
\medskip

\noindent In other words, ${\sf rank}_S(x)$ is the {\em regular 
$\mathcal J$-depth} of the {\em idempotent generated by} $x$.

We also define the rank with respect to automata, i.e., we extend
${\sf rank}_S: S \to {\mathbb N}$ to a {\em word-rank} function 
${\sf rank}_{\bf A}: A^* \to {\mathbb N}$ \, (where ${\mathbb N} = $
$\{0, 1, \ldots \}$ denotes the natural integers). 
For an automaton {\bf A} over $A$ and $w \in A^*$, let $[w]$ be the image 
of $w$ in the syntactic monoid $S$ of {\bf A}. Then ${\sf rank}_{\bf A}(w)$ 
is defined to be ${\sf rank}_S([w])$.

\begin{defn} {\bf (Rank condition).}  \label{defRankCond} \ 
An automaton morphism $\varphi: {\bf A} \twoheadrightarrow {\bf B}$ in 
${\mathcal A}_A$ satisfies the {\em rank condition} iff the following holds.
For all $q \in Q_{\bf A}$ and all $\alpha, \beta \in A^+$ such that
$q \cdot_{_{\bf A}} \alpha = q$ and 
$(q)\varphi \cdot_{_{\bf B}} \beta = (q)\varphi$ and 
${\sf rank}_{\bf B}(\alpha) \geq {\sf rank}_{\bf B}(\beta)$, 
we have: \ $q \cdot_{_{\bf A}} \beta = q$. 
\end{defn}
The formal statement of the rank condition on $\varphi$ is as follows.

\medskip

$(\forall q \in Q_{\bf A}) (\forall \alpha, \beta \in A^+):$

 \hspace{.5in}  $[ \, q \cdot_{_{\bf A}} \alpha = q$ \ \ $\&$ 
 \ \ $(q)\varphi \cdot_{_{\bf B}} \beta = (q)\varphi$ \ \ $\&$
 \ \ ${\sf rank}_{\bf B}(\alpha) \geq {\sf rank}_{\bf B}(\beta) \, ]$
 \ \ implies \ \ $q \cdot_{_{\bf A}} \beta = q$.

\medskip

\noindent
Yet another way to say this: If $\alpha$ fixes $q$ in {\bf A}, and $\beta$ 
fixes $(q)\varphi$ in $({\bf A})\varphi$, and 
${\sf rank}_{({\bf A})\varphi}(\alpha)$ $\geq$
${\sf rank}_{({\bf A})\varphi}(\beta)$, then 
$\beta$ also fixes $q$ in {\bf A}. 

\medskip

\noindent {\bf The following is a better statement of the rank condition:}

\medskip

$(\forall q \in Q_{\bf A}) (\forall \beta \in A^+):$
 
 if \ \ $[ \, (q)\varphi \cdot_{_{\bf B}} \beta = (q)\varphi$ \ \ $\&$ \ \   
 $(\exists \alpha \in A^+)$
    $[ \, q \cdot_{_{\bf A}} \alpha = q$ \ $\&$ \   
    ${\sf rank}_{\bf B}(\alpha) \geq {\sf rank}_{\bf B}(\beta) \, ] \, ]$
  \ \ then \ \ $q \cdot_{_{\bf A}} \beta = q$.

\medskip

\noindent In words: \ If $\beta$ fixes $(q)\varphi$ in $({\bf A})\varphi$, 
and there is $\alpha$ that fixes $q$ in {\bf A} such that 
${\sf rank}_{({\bf A})\varphi}(\alpha)$ $\geq$
${\sf rank}_{({\bf A})\varphi}(\beta)$, then
$\beta$ also fixes $q$ in {\bf A}.

\begin{lem}
If two morphisms $\varphi_1: {\bf A}_1 \to {\bf B}$ and
$\varphi_2: {\bf A}_2 \to {\bf B}$ satisfy the rank condition, 
then the direct product morphism
$\varphi_1 \times_A \varphi_2: {\bf A}_1 \times_A {\bf A}_2 \to {\bf B}$
in ${\mathcal A}_A$ also satisfies the rank condition.
\end{lem}
{\bf Proof.} Assume that \    
$(I_1 \cdot_1 w, \ I_2 \cdot_2 w) \cdot \alpha$ $=$
$(I_1 \cdot_1 w, \ I_2 \cdot_2 w)$ \ in 
${\bf A}_1 \times_A {\bf A}_2$, assume that
 \ $I \cdot_{_{\bf B}} w \cdot_{_{\bf B}} \beta = I \cdot_{_{\bf B}} w$
 \ in {\bf B}, and assume that 
${\sf rank}_{\bf B}(\alpha) \geq {\sf rank}_{\bf B}(\beta)$.
The equality in ${\bf A}_1 \times_A {\bf A}_2$ implies that for $i = 1,2:$
 \ $I_i \cdot_i w \cdot_i \alpha = I_i \cdot_i w$ and 
$\varphi_i(I_i) \cdot_{_{\bf B}} w \cdot_{_{\bf B}} \beta$ $=$
$\varphi_i(I_i) \cdot_{_{\bf B}} w$.  Hence by the rank condition for
$\varphi_i$ we have $I_i \cdot_i w \cdot_i \beta = I_i \cdot_i w$. So,
$\beta$ fixes $(I_1 \cdot_1 w, \ I_2 \cdot_2 w)$ in 
${\bf A}_1 \times_A {\bf A}_2$. 
 \ \ \ $\Box$

\subsection{Definition of the glc-cover}

Let {\sf cov} be any cover in the category ${\mathcal A}_A$.
For an automaton {\bf A} in ${\mathcal A}_A$ we consider the interval
$[{\bf A}, {\bf A}^{\sf cov}]$. This is a finite lattice, as we saw. 

\begin{defn}
The {\em glc}-cover ${\bf A}^{\sf cov \, glc}$ is defined by the following
join in the category ${\cal A}_A$: 

\medskip

${\bf A}^{\sf cov \, glc} \ = \ $
$\bigvee \, \{ {\bf X} \in [{\bf A}, {\bf A}^{\sf cov}] : $ there exists a 
morphism ${\bf A} \twoheadleftarrow {\bf X}$ in ${\mathcal A}_A$ that is 
1:1${\mathcal R}$ and that  

\hspace{1.98in}  satisfies the rank condition\}.
\end{defn}
More generally, if ${\bf A}$ and ${\bf B}$ are $A$-automata such that
${\bf A} \twoheadleftarrow {\bf B}$, we can define the {\em glc}-cover
of the interval $[{\bf A}, {\bf B}]$ by

\medskip

$[{\bf A}, {\bf B}]^{\sf glc} \ = \ $
$\bigvee \, \{ {\bf X} \in [{\bf A}, {\bf B}] : $ 
{\it there exists a morphism ${\bf A} \twoheadleftarrow {\bf X}$ in 
    ${\mathcal A}_A$ that is 1:1${\mathcal R}$ and that

\hspace{1.9in}  satisfies the rank condition\}.
}  

\medskip

\noindent By the previous Lemmas, and since $[{\bf A}, {\bf A}^{\sf cov}]$ 
and $[{\bf A}, {\bf B}]$ are finite complete lattices, these joins exist.  
The name ``glc'' stands for global-local-cover.  

Sometimes we use the glc-cover without the rank condition, i.e., with 
the 1:1${\mathcal R}$ requirement alone; 
we will make that clear in the context. 

\medskip

When the automata {\bf A} and {\bf B} in $[{\bf A}, {\bf B}]^{\sf glc}$ 
are right Cayley graphs it makes sense to define the cover
$[{\bf A}, {\bf B}]^{\sf glc}$ either as above (letting {\bf X} range over
$A$-automata), or just in terms of right Cayley graphs (i.e., letting 
{\bf X} range over right Cayley graphs only).  We will prove next that the 
resulting cover is the same in either approach.

Let {\bf A} and {\bf B} be right Cayley graphs of $A$-generated semigroups, 
and suppose {\bf C} is any $A$-automaton such that 
${\bf B} \twoheadrightarrow {\bf C}$
$ \stackrel{\varphi}{\twoheadrightarrow} {\bf A}$, where $\varphi$ is 
1:1${\mathcal R}$ (with or without the rank condition). Let $S({\bf C})$ 
be the $A$-generated syntactic semigroup of {\bf C}, and let 
$\Gamma(S({\bf C})) = \Gamma_A(S({\bf C}),\sigma)$ be the right Cayley graph 
of $S({\bf C})$, where $\sigma$ embeds $A$ into a generating set of 
$S({\bf C})$.  Let $\eta: \Gamma(S({\bf C})) \twoheadrightarrow {\bf C}$ be 
the canonical surmorphism, mapping semigroup elements to states.

\begin{lem}  \ If \ ${\bf B} \twoheadrightarrow \Gamma(S({\bf C}))$
$\stackrel{\eta}{\twoheadrightarrow} {\bf C}$
$\stackrel{\varphi}{\twoheadrightarrow} {\bf A}$, where {\bf A} and 
{\bf B} are Cayley graphs, and if $\varphi$ is 1:1${\mathcal R}$, then 
$\eta \circ \varphi: \Gamma(S({\bf C})) \twoheadrightarrow {\bf A}$ is 
also 1:1${\mathcal R}$.
If, in addition, $\varphi$ satisfies the rank condition then 
$\eta \circ \varphi$ also satisfies the rank condition.  
\end{lem}
{\bf Proof.} If for $s_1, s_2 \in S({\bf C})$ we have $s \ne s_2$, 
$s_1 \equiv_{\mathcal R} s_2$, then by faithful action of $S({\bf C})$ 
in {\bf C} there exists a state $q$ of {\bf C} with 
$q \cdot s_1 \neq q \cdot s_2$. 
Hence since $\varphi$ is 1:1${\mathcal R}$, 
$(q s_1)\varphi \ne (q s_2)\varphi$.
Hence, since {\bf A} is a Cayley graph, 
$I_{\bf A} \cdot s_1 \ne I_{\bf A} \cdot s_2$, so 
$(s_1)\eta \circ \varphi$ $ = $ $I_{\bf A} \cdot s_1$ $\ne$ 
$I_{\bf A} \cdot s_2 = (s_2)\eta \circ \varphi$, i.e.,
$\eta \circ \varphi$ is 1:1${\mathcal R}$.  

Suppose $\varphi$ satisfies the rank condition. We want to show for all
$s \in S({\bf C})$ and $\alpha, \beta \in A^+ :$ \ if
 \, $s \cdot_{_S} \alpha = s$  and 
 \ $(s)\eta \varphi \cdot_{_{\bf A}} \beta = (s)\eta \varphi$ and 
${\sf rank}_{\bf B}(\alpha) \geq {\sf rank}_{\bf B}(\beta)$, 
then $s \cdot_{_S} \beta = s$. Here $\cdot_{_S}$ is the next-state function
in $\Gamma(S({\bf C}))$.
By the rank condition for $\varphi$, we have 
 \ $(s)\eta \cdot_{_{\bf C}} \beta = (s)\eta$. Hence, since $\eta$ is 
bijective on the semigroups,  $s \cdot_{_{\bf C}} \beta = s$.
 \ \ \ $\Box$

\begin{pro} \label{glcAutomVsCayley} \ 
When {\bf A} and {\bf B} are $A$-generated Cayley graphs then the cover
$[{\bf A}, {\bf B}]^{\sf glc}$ is the same in the category of $A$-automata
as in the category of $A$-generated Cayley graphs.
\end{pro}
{\bf Proof.} Since an $A$-generated Cayley graph is an $A$-automaton,
the $A$-automata cover (which is the join over automata) maps onto the
Cayley graph cover.

By the previous lemma, the Cayley graph of the automaton cover maps
onto {\bf A} by a 1:1${\mathcal R}$ morphism with rank condition. Hence 
the Cayley graph is equal to the cover, by maximality of the cover.  
 \ \ \ $\Box$


\subsection{Examples of {\sf glc}-covers for 1:1${\mathcal R}$ morphisms}

\noindent {\bf (1)} \ Let {\bf V} be a pseudovariety of semigroups, 
and let $S_A$ be an $A$-generated semigroup. Let $S_A^{\bf V}$ denote
the {\it maximum image of $S_A$ in} {\bf V}; i.e., 

\smallskip

 \ \ \  $S_A^{\bf V} \ = \ $
$\bigvee \, \{ Y \in {\bf V} : $ there exists a morphism 
$S_A \twoheadrightarrow Y$ in ${\mathcal S}_A\}$, 

\smallskip

\noindent where $\vee$ is the join in the category of $A$-generated 
semigroups; this join is finite since $S_A$ is finite.   

Let {\bf R} be the pseudovariety of $\mathcal R$-trivial semigroups. 
For any $A$-generated semigroup $S_A$, $S_A^{\bf R}$ denotes thus the 
maximum $\mathcal R$-trivial image. References for this and the next 
few definitions are \cite{qTheory} and \cite{GST}. 

Let \mal be the {\it Maltsev product}. The {\it Maltsev kernel} of a 
morphism $\psi: S \twoheadrightarrow T$ is, by definition, the 
pseudovariety generated by the set of semigroups  
$\{ (e)\psi^{-1} : e = e^2 \in T\}$.
We consider only the case where {\bf V} is {\it locally finite} (i.e., every 
finitely generated semigroup in {\bf V} is finite). For example, the 
pseudo-variety $\langle T \rangle$ generated by any finite semigroup $T$ is locally finite.

Then {\bf V} \mal ({\bf .}) is an expansion of $A$-generated semigroups, 
where {\bf V} \mal $S_A$ is the unique functorially largest 
$A$-generated semigroup that maps onto $S_A$ by a surmorphism whose 
Maltsev kernel is in {\bf V}: 

\medskip

 \ \ \  {\bf V} \mal $S_A$ $ \ = \ $   
$\bigvee \, \{ X : $ there exists $\varphi: X \twoheadrightarrow S_A$ in 
${\cal S}_A$ such that the Maltsev kernel of $\varphi$ is  

 \hspace{1.5in} contained in {\bf V}\}.  

\medskip

\noindent Since {\bf V} is locally finite, this join is a finite set, so 
{\bf V} \mal $S_A$ is finite (by Brown's theorem; see \cite{GST},
\cite{qTheory}).

\begin{pro} \label{examplesGLC} {\bf (Example \ref{examplesGLC})}  
 \ Let $\Gamma(S_A) = [\Gamma(T_A), \Gamma(U_A)]^{\sf glc}$ be the 
{\sf glc}-cover of $A$-generated right Cayley graphs with respect to the 
1:1${\mathcal R}$ property; here we temporarily drop the rank condition. 
Hence, $S_A$ is the unique maximum $A$-generated semigroup such that

\smallskip

\hspace{1in} 
$\Gamma(U_A) \stackrel{\varphi_1}{\twoheadrightarrow} \Gamma(S_A)$
$\stackrel{\varphi_2}{\twoheadrightarrow} \Gamma(T_A)$, 

\smallskip

\noindent where $\varphi_2$ is 1:1${\mathcal R}$.  Then: 

\smallskip

\noindent {\bf (a)} \ The Maltsev kernel of $\varphi_2$ (and in fact, of
any 1:1${\mathcal R}$ morphism) is contained in {\bf R}.

\smallskip

\noindent {\bf (b)} \ If $T_A \in {\bf R}$ and $U_A = T_A^{\sf cov}$
for some cover \, $(.)^{\sf cov}$, then

\smallskip

\hspace{1in}
$[\Gamma(T_A), \Gamma(T_A^{\sf cov})]^{\sf glc}$ $ \ = \ $
$T_A^{ {\sf cov} \, {\bf R}}$.

\smallskip

\noindent {\bf (c)} \ For any $A$-generated semigroup $W_A$ we consider 
the expansion $(.)_A^{\sf exp} = \langle W_A \rangle$ \mal $({\bf .})$,
where $\langle W_A \rangle$ is the pseudovariety generated by $W_A$. 
Let ${\bf 1}_A$ be the one-element semigroup in the category of
$A$-generated semigroups.

\smallskip

Then ${\bf 1}_A^{\sf exp} = W_A$, and
 \, $[\Gamma({\bf 1}_A), \Gamma(W_A)]^{\sf glc} = W_A^{\bf R}$. 

Thus, for any $A$-generated  $R \in {\bf R}$, $\langle R \rangle$ can be 
the Maltsev kernel.
\end{pro}
{\bf Proof.} {\bf (a)} \ If \     
$\varphi_2: \Gamma(S_A) \twoheadrightarrow \Gamma(T_A)$ \ is
1:1${\mathcal R}$ then the corresponding map $S_A \twoheadrightarrow T_A$
(which we also denote by $\varphi_2$) is also 1:1${\mathcal R}$. 
Then for any idempotent $e = e^2 \in T_A$, $(e)\varphi_2^{-1}$ is an
${\mathcal R}$-trivial semigroup. Indeed, if $x_1 \equiv_{\mathcal R} x_2$ 
in $(e)\varphi_2^{-1}$ then $(x_1)\varphi_2 = e = (x_2)\varphi_2$, hence 
the 1:1${\mathcal R}$-property implies $x_1 = x_2$. 

{\bf (b)} \ When
$\varphi_2: S_A \twoheadrightarrow T_A$ is 1:1${\mathcal R}$, then 
$S_A \in {\bf R}$ iff $T_A \in {\bf R}$. 
Indeed, if $x_1 \equiv_{\mathcal R} x_2$ in $S_A$ then
$(x_1)\varphi_2 \equiv_{\mathcal R} (x_2)\varphi_2$ in $T_A$, hence 
$(x_1)\varphi_2 = (x_2)\varphi_2$ (since $T_A \in {\bf R}$), hence
$x_1 = x_2$ (since $\varphi_2$ is 1:1${\mathcal R}$);
the converse is obvious.

Since $[\Gamma(T_A), \Gamma(T_A^{\sf  exp})]^{\sf glc}$
$\twoheadrightarrow \Gamma(T_A)$ is 1:1${\mathcal R}$ and since
$T_A \in {\bf R}$, it follows that 
$[\Gamma(T_A), \Gamma(T_A^{\sf  exp})]^{\sf glc}$ is in {\bf R}.
Since the {\sf glc}-cover is maximal with respect to the
1:1${\mathcal R}$-property, it follows that the {\sf glc}-cover is 
$T_A^{ {\sf  exp} \, {\bf R}}$ (by the definition of $(.)^{\bf R}$). 
 
{\bf (c)} \ Obviously, the map $S_A \twoheadrightarrow {\bf 1}_A$ is 
1:1${\mathcal R}$ iff $S_A \in {\bf R}$.
Hence, by the definition of the expansion $(.)^{\sf exp}$ $=$
$\langle W_A \rangle$ \mal $({\bf .})$ we have 
${\bf 1}_A^{\sf exp} = W_A$. And by the definition of $(.)^{\bf R}$, 
 \ $[\Gamma({\bf 1}_A), \Gamma(W_A)]^{\sf glc} = W_A^{\bf R}$.     

\smallskip 

Now let $W_A$ ($ = R$) be any $A$-generated semigroup in {\bf R}. Then 
$W_A = W_A^{\bf R}$, so in that case, \ 
$[\Gamma({\bf 1}_A), \Gamma(W_A)]^{\sf glc} = W_A$ \ (by {\bf (c)}).

And the Maltsev kernel is $\langle W_A \rangle$ (which is 
$\langle R \rangle$).
 \ \ \ $\Box$

\medskip

\noindent {\bf Remark.} From Example \ref{examplesGLC} we see that 
$[\Gamma_A(T), \Gamma_A(T^{\sf cov})]^{\sf glc}$ can be:
 
(case 1) \ $\Gamma_A(T)$, or 
 
(case 2) \ $\Gamma_A(T^{\sf cov})$, or 

(case 3) \ strictly between the bounds of the interval.

\smallskip

\noindent Let $T_A = {\bf 1}_A$ and 
$(.)_A^{\sf cov} = \langle W_A \rangle$ \mal $({\bf .})$.  Then 
${\bf 1}_A^{\sf cov} = W_A$, and by Example \ref{examplesGLC}.(c), 
$[\Gamma_A(T), \Gamma_A(T^{\sf cov})]^{\sf glc}$ $=$ 
$[{\bf 1}_A, {\bf 1}_A^{\sf cov}]^{\sf glc} = W_A^{\bf R}$. 
 
\smallskip

\noindent (1) When $W_A$ is a non-trivial $A$-generated group, then 
$W_A^{\bf R} = {\bf 1}_A$. So we get (case 1). 

\smallskip

\noindent (2) When ${\bf 1}_A \neq W_A \in {\bf R}$ (e.g., $W_A$ is the 
left-zero semigroup generated by $A$) then by the Example 
\ref{examplesGLC}.(b), we are in (case 2) and not at the same time in 
(case 1). 

\smallskip

\noindent (3) Let $W_A \not\in {\bf R}$ and $W_A^{\bf R} \neq {\bf 1}$; e.g.,
let $W_A = {\mathbb Z}_2 \times_A A^{\ell}$, where $A^{\ell}$ is the left-zero
semigroup generated by $A$. Then $W_A \twoheadrightarrow A^{\ell}$ by
projection; then $[{\bf 1}_A, {\bf 1}_A^{\sf cov}]^{\sf glc} = W_A^{\bf R}$ 
is strictly between ${\bf 1}_A$ and ${\bf 1}_A^{\sf cov} = W_A$.     
This illustrates (case 3).

\bigskip

\noindent {\bf (2)} \ Let us consider the reverse of the {\it delay} 
pseudovariety ${\bf D} = [x y^{\omega} = y^{\omega}]$. 
So ${\bf D}^{\rm rev} = [ x^{\omega} y = x^{\omega}]$, corresponding to
the limit of the set of equations  
$\{ x_1 \ldots x_k x_{k+1} = x_1 \ldots x_k : k \geq 1\}$.
It is known that $S \in {\bf D}^{\rm rev}$ iff $S$ is a nilpotent extension
of a left-zero semigroup; such a semigroup consist of a minimal ideal which
is an $\cal L$-class $L$ of left-zeros, and every element of the semigroup
satisfies $(\exists k) \ x^k \in L$. The left-zero semigroups form the 
pseudovariety {\bf LZ} $=$ $[x y = x]$. For fixed $k > 0$, let 
${\bf D}_k^{\rm rev} = [x_1 \ldots x_k x_{k+1} = x_1 \ldots x_k]$ be the 
reverse of the delay-$k$ pseudovariety; this pseudovariety is locally finite.
We consider the expansion $(.)^{\sf exp} = {\bf D}_k^{\rm rev} \mal (.)$.
Then we have for every $A$-generated semigroup $T_A$: 

\smallskip

 \ \ \ {\it  ${\bf D}_k^{\rm rev} \mal T_A$ \, $=$ \, 
$T_A^{\sf exp} \ \twoheadrightarrow \ T_A$ \ is injective on every 
regular ${\mathcal R}$-class. 
}

\smallskip

\noindent Indeed, let us apply Rees' Theorem to a regular 
${\mathcal R}$-class $\{ a \} \times G \times B$ of $T_A^{\sf exp}$. If the 
map identifies two elements in this ${\mathcal R}$-class then it identifies 
two $\cal L$-classes, hence it identifies to ${\mathcal L}$-equivalent
idempotents. So in that case there are two idempotents $e \ne f$. Thus, 
$\{e, f\}$ is a right-zero semigroup, but a non-trivial right-zero 
semigroup cannot be in ${\bf D}_k^{\rm rev}$. 

\medskip

\noindent {\bf Notation:} A morphism that is injective on every regular 
$\cal R$-class is said to be 1:1reg$\cal R$.

\bigskip 

\noindent {\bf Digression: \, Example of a morphism 
$\varphi: S_A \twoheadrightarrow T_A$ that is 
1:1{\rm reg}$\cal R$, but not 1:1$\cal R$ } 

\smallskip

Let $S = \{\i, x, \ov{\i}, \ov{x}, 0\}$, \ $T = \{\i, x, y, 0\}$, 
 \ $\varphi = \{(\i,\i), (x,x), (\ov{\i}, y), (\ov{x}, y), (0,0)\}$, and
$A$ is a two-element set embedded injectively into $\{x, \ov{\i}\}$ and
$\{x, y\}$.  The multiplication tables of $S$ and $T$ are: 

\bigskip

\hspace{.7in}
\begin{tabular}{r||r|r|r|r|r}
$\bullet_{_S}$ & $\i$ & $x$   & $\ov{\i}$ & $\ov{x}$ & 0 \\ \hline \hline 
$\i$          & $\i$ & $x$    & $\ov{\i}$ & $\ov{x}$ & 0 \\ \hline
$x$           &  $x$ & $\i$   & $\ov{\i}$ & $\ov{x}$ & 0 \\ \hline
$\ov{\i}$ & $\ov{\i}$ & $\ov{x}$& 0       & 0        & 0 \\ \hline
$\ov{x}$  & $\ov{x}$  & $\ov{\i}$& 0      & 0        & 0 \\ \hline
0         & 0         & 0     & 0         & 0        & 0 \\ 
\end{tabular}
\hspace{1in}
\begin{tabular}{r||r|r|r|r}
$\bullet_{_T}$ & $\i$ & $x$  & $y$  & 0  \\ \hline \hline    
$\i$           & $\i$ & $x$  & $y$  & 0  \\ \hline
$x$            & $x$  & $\i$ & $y$  & 0  \\ \hline
$y$            & $\i$ & $y$  & 0    & 0  \\ \hline
0              & 0    & 0    & 0    & 0  \\ 
\end{tabular}

\bigskip

\noindent The structure of these semigroups can be understood as follows:
$S$ has three $\cal J$-classes, namely the group 
${\mathbb Z}_2 = \{\i, x\}$ at the top, a null $\cal R$-class 
$\{ \ov{\i}, \ov{x}\}$ in the middle, and a zero 0. The group  
${\mathbb Z}_2$ acts on the $\cal R$-class in the obvious way on the left 
and on the right (see the multiplication table), and ``null'' means that 
products in the middle $\cal R$-class are 0.  The semigroup 
$T$ differs from $S$ only in that its middle null $\cal J$-class is a
singleton.  The action of ${\mathbb Z}_2$ and the null property of the 
middle $\cal J$-class make it easy to check associativity. 

It is clear that $\varphi$ is injective on regular $\cal R$-classes (namely
${\mathbb Z}_2$ and $\{0\}$), and not injective on the null middle
$\cal R$-class.
 \ \ \ [End, digression] 

\bigskip

We continue to use {\sf glc}-covers based on 1:1$\cal R$ morphisms,
ignoring the rank condition. We also consider 1:1reg$\cal R$ morphisms, and
we distinguish the two by the notation $(.)^{{\sf glc} \, 1:1{\cal R}}$,
respectively $(.)^{{\sf glc} \, 1:1{\rm reg}{\cal R}}$.  
We will use the following.

\begin{lem}
Let $\alpha, \beta$ be surmorphisms of $A$-generated semigroups. Then 
$\alpha \circ \beta$ is {\rm 1:1}${\cal R}$ (or {\rm 1:1reg}$\cal R$)
iff each of $\alpha$ and $\beta$ is {\rm 1:1}${\cal R}$ (respectively 
{\rm 1:1reg}$\cal R$).
\end{lem}
{\bf Proof.} The right-to-left implication is obvious. Assume now that
$\beta \circ \alpha$ is 1:1$\cal R$, from which it immediately follows that
$\alpha$ is 1:1$\cal R$. For $\beta$, let $R$ be an $\cal R$-class of
${\sf dom}(\beta)$, and let $y_1, y_2 \in R$ with $y_1 \ne y_2$. Let
$x_1, x_2 \in {\sf dom}(\alpha)$ be such that $x_1 \equiv_{\cal R} x_2$
and $y_1 = (x_1)\alpha$ and $y_2 = (x_2)\alpha$; $x_1, x_2$ exist since
$\alpha$ is surjective, and since the inverse image of
an $\cal R$-class is a union of $\cal R$-classes. Then
$(y_1)\beta = (x_1)\alpha \beta \ne (x_2)\alpha \beta = (y_2)\beta$,
where ``$\ne$'' holds because $(.)\alpha \beta$ is 1:1$\cal R$.
Hence $(y_1)\beta \ne (y_2)\beta$, i.e., $\beta$ is 1:1$\cal R$.
The same proof (when $R$ is regular) works for 1:1reg$\cal R$, using
the fact that the inverse image of a regular $\cal R$-class contains a
regular $\cal R$-class.
 \ \ \ $\Box$

\bigskip

\noindent We consider the expansion defined by \, $T_A^{\sf exp}$   
$ \ = \ $ ${\bf D}_k^{\rm rev}$ \mal $T_A \ \ \twoheadrightarrow \ \ T_A$. 
Then  

\smallskip

 \ \ \ $T_A^{\sf exp} \ \ \twoheadrightarrow \ \  $ 
$[T_A, T_A^{\sf exp}]^{{\sf glc} \ 1:1{\rm reg}{\cal R}}$

\smallskip
 
\noindent is an isomorphism, since $T_A^{\sf exp} \twoheadrightarrow T_A$
is 1:1reg$\cal R$. So we have 

\smallskip

 \ \ \  $T_A^{\sf exp}$ $ \ = \ $
$[T_A, T_A^{\sf exp}]^{{\sf glc} \ 1:1{\rm reg}{\cal R}}$
$ \ \ \ \stackrel{\varphi_1}{\twoheadrightarrow} \ \ \ $
$[T_A, T_A^{\sf exp}]^{{\sf glc} \ 1:1{\cal R}}$
$ \ \ \ \stackrel{\varphi_2}{\twoheadrightarrow} \ \ \ T_A$, 

\smallskip  

\noindent where $\varphi_1$ is not an isomorphism (in general), because 
1:1reg$\cal R$ does not imply 1:1$\cal R$.
The following is straightforward.

\begin{pro}
Let $(.)^{\sf E}$ be an expansion such that 
$S_A^{\sf \ E} \twoheadrightarrow S_A$ is {\rm 1:1}$\cal R$ for all $S_A$.
If the semigroups $U_A$ and $T_A$ satisfy $U_A^{\sf \ E} = U_A$ and 
$T_A^{\sf \ E} = T_A$, then 
 \, $[T_A, U_A]^{ {\sf glc} \ {\sf E}} = [T_A, U_A]^{\sf glc}$.  
 \ \ \ $\Box$
\end{pro}
For example, the right Rhodes expansion $(.)^{\wedge R}$ has the 1:1$\cal R$ 
property used above.

\subsection{The ${\mathcal R}$- and ${\mathcal L}$-expansions of an 
automaton}

\noindent {\bf The ${\mathcal R}$-expansion:} \   We generalize the right 
Rhodes ${\wedge_{\mathcal R}}$-expansion to the automaton category 
${\mathcal A}_A$. In the special case of the right Cayley graph automaton 
$(S^I, I, \cdot)$ of a semigroup $S_A$ over the alphabet $A$ we recover the 
known Rhodes expansion $S_A^{\wedge_{\mathcal R}}$.  

The definition proceeds in two steps. 
Let ${\bf A} = (Q, i, \cdot)$ be an automaton in ${\mathcal A}_A$ with 
start state $i$. 
We assume that the start state $i$ is not reachable from any other state.

\noindent (1) Let $Q^{\wedge_{\mathcal R}} = \ $ 
$\{ (i >_{\mathcal R} q_1 >_{\mathcal R} \ \ldots \ >_{\mathcal R} q_k) : $
$k \geq 0 \ {\rm and} \ q_1, \ldots, q_k \in Q\}$ \ be the set of strict 
reachability chains, starting at the start state.  From this we build an 
automaton with state set $Q^{\wedge_{\mathcal R}}$, 
start state $(i)$, and next-state function $\bullet$; the latter is defined 
as follows. For any $a \in A$ and ${\bf q} = $
$(i >_{\mathcal R} q_1 >_{\mathcal R} \ \ldots \ >_{\mathcal R} q_{k-1}$
$>_{\mathcal R} q_k)$ $\in Q^{\wedge_{\mathcal R}}$,  

\smallskip
 
 \ \ \ ${\bf q} \bullet a \ = \ $
$(i >_{\mathcal R} q_1 >_{\mathcal R} \ \ldots \ >_{\mathcal R} q_{k-1} $
$ >_{\mathcal R} q_k >_{\mathcal R} q_k \cdot a)$, 
 \ \ if \ $q_k >_{\mathcal R} q_k \cdot a$; \ \ and

\smallskip

 \ \ \ ${\bf q} \bullet a \ = \ $
$(i >_{\mathcal R} q_1 >_{\mathcal R} \ \ldots \ >_{\mathcal R} q_{k-1}$
$>_{\mathcal R} q_k \cdot a)$, \ \ if 
 \ $q_k \equiv_{\mathcal R} q_k \cdot a$.

\medskip

\noindent 
(2) Let $Q_A^{\wedge_{\mathcal R}} \subseteq Q^{\wedge_{\mathcal R}}$ be 
the set of states of the above automaton that are reachable from the start
state $(i)$; i.e., 
$Q_A^{\wedge_{\mathcal R}} = \{ (i) \bullet w \ : \ w \in A^*\}$.   
We define ${\bf A}_A^{\wedge_{\mathcal R}}$ to be the sub-automaton 
$(Q_A^{\wedge_{\mathcal R}}, \ (i), \ \bullet)$ of the above automaton.

The canonical natural transformation $\eta$ of the expansion is defined for 
each {\bf A} by $\eta_{\bf A}: {\bf A}_A^{\wedge_{\mathcal R}} \to {\bf A}$, 
where \ $(i >_{\mathcal R} q_1 >_{\mathcal R} \ \ldots \ $
$>_{\mathcal R} q_k)\eta_{\bf A} \ = \ q_k$. 
One easily checks that this is an automaton morphism in the category 
${\mathcal A}_A$. 

\smallskip

Any automaton morphism $\varphi: {\bf A} \to {\bf B}$ in ${\mathcal A}_A$
can be expanded to 
$\varphi_A^{\wedge_{\mathcal R}}:$
${\bf A}_A^{\wedge_{\mathcal R}} \to {\bf B}_A^{\wedge_{\mathcal R}}$ 
defined by 

\smallskip

 \ \ \ \   
$(i_{\bf A} >_{\mathcal R} q_1  >_{\mathcal R} \ \ldots \ $
$ >_{\mathcal R} q_k)\varphi_A^{\wedge_{\mathcal R}}$ \ $=$ \ 
${\sf red}\big( i_{\bf B} >_{\mathcal R} (q_1)\varphi$
$\geq_{\mathcal R} \ \ldots \ \geq_{\mathcal R} (q_k)\varphi \big)$.
 
\smallskip

\noindent Here, the ``reduction operation'' {\sf red} has the effect of
replacing every maximal $\equiv_{\mathcal R}$-chain by its rightmost 
element. More formally, ${\sf red}({\bf q}) = {\bf q}$ if {\bf q} is a 
strict chain; and \ 

\smallskip

${\sf red}( \ \ldots \ \geq_{\mathcal R} q_{i-1} \geq_{\mathcal R} $
$q_i \equiv_{\mathcal R} q_{i+1} \geq_{\mathcal R} \ \ldots \ )$ \ $=$ 
 \ ${\sf red}( \ \ldots \ \geq_{\mathcal R} q_{i-1} \geq_{\mathcal R} $
$q_{i+1} \geq_{\mathcal R} \ \ldots \ )$.

\smallskip

\noindent It is easy to check that $\varphi_A^{\wedge_{\mathcal R}}$ 
is an automaton morphism. 

In summary, $(.)_A^{\wedge_{\mathcal R}}$ is an expansion in the category 
of automata ${\mathcal A}_A$, according to the categorial definition of 
expansions. 

\smallskip

The analogy between the Rhodes $\wedge_{\mathcal R}$-expansion and the 
automata $\wedge_{\mathcal R}$-expansion goes beyond the similarity of 
the definitions.  
Let $S_A$ be the syntactic monoid of an $A$-automaton 
${\bf A} = (Q, i, \cdot)$; hence $(Q, S_A)$ is a faithful right action. 

\begin{pro}
The Rhodes expansion $S_A^{\wedge_{\mathcal R}}$ acts on the set
$Q_A^{\wedge_{\mathcal R}}$ by the following right action.
\end{pro} 
For ${\bf q} =$
$(i >_{\mathcal R} q_1 >_{\mathcal R} \ \ldots \ >_{\mathcal R} q_k)$
$ \in Q_A^{\wedge_{\mathcal R}}$ and ${\bf s} =$ 
$(1 >_{\mathcal R} s_1 >_{\mathcal R} >_{\mathcal R} s_2 >_{\mathcal R}$
$ \ \ldots \ >_{\mathcal R} s_h) \in S_A^{\wedge_{\mathcal R}}$ this 
action is defined by:

\smallskip

${\bf q} \bullet {\bf s} \ = \  $ 
${\sf red}(i >_{\mathcal R} q_1 >_{\mathcal R} \ \ldots \ >_{\mathcal R} $
$q_k \geq_{\mathcal R} q_k \cdot s_1 \geq_{\mathcal R} q_k \cdot s_2 $
$ \geq_{\mathcal R}  \ \ldots \ \geq_{\mathcal R} q_k \cdot s_h)$.  
 
\medskip

\noindent This action yields a homomorphism from $S_A^{\wedge_{\mathcal R}}$ 
onto the syntactic monoid of ${\bf A}_A^{\wedge_{\mathcal R}}$. 
The action is not necessarily faithful, i.e., this homomorphism is not 
necessarily an isomorphism. 
Just as for the Rhodes expansion for the semigroup category ${\mathcal S}_A$,
we have in the automaton category ${\mathcal A}_A$:
\begin{pro}
 \ For every automaton {\bf A} in ${\mathcal A}_A:$ \
$({\bf A}_A^{\wedge_{\mathcal R}})_A^{\wedge_{\mathcal R}} $
$= {\bf A}_A^{\wedge_{\mathcal R}}$.
\end{pro}
{\bf Proof.} The proof is similar to the case of
$S_A^{\wedge_{\mathcal R}}$.
 \ \ \ $\Box$

\begin{pro}
The reachability order among the states of 
${\bf A}_A^{\wedge_{\mathcal R}}$ is unambiguous.  
\end{pro}
{\bf Proof.} The proof is similar to the case of
$S_A^{\wedge_{\mathcal R}}$.
 \ \ \ $\Box$

\begin{pro}
For any automaton {\bf A}, the canonical map 
$\eta_{\bf A}: {\bf A}_A^{\wedge_{\mathcal R}} \to {\bf A}$ 
is fully 1:1${\mathcal R}$ and satisfies the rank condition.
\end{pro}
{\bf Proof.} \ [1:1${\mathcal R}$] \ If ${\bf p} \equiv {\bf q} \in$
$A^{\wedge_{\mathcal R}}$ then {\bf p} and {\bf q} have the form
${\bf p} = (r_1 > \ldots > r_{k-1} > p)$, respectively
${\bf q} = (r_1 > \ldots > r_{k-1} > q)$. Hence, when 
${\bf p} \equiv {\bf q}$ we have: \ ${\bf p} \neq {\bf q}$ \ iff 
 \ $\eta({\bf p}) \neq \eta({\bf q})$.

\smallskip

{\sf [fully 1:1${\mathcal R}$]} \ Let $p = \eta({\bf p})$ and suppose 
$p \equiv q$ in {\bf A}. We need to show that there is {\bf q} such that 
$\eta({\bf q}) = q$ and ${\bf p} \equiv {\bf q}$. 
Since $p = \eta({\bf p})$, the chain {\bf p} has the form 
$(r_1 > \ldots > r_{k-1} > p)$ for some $r_1, \ldots, r_{k-1} \in Q$.
Since $p \equiv q$, the chain $(r_1 > \ldots > r_{k-1} > q)$ (let's call it
{\bf q}) exists in $Q^{\wedge_{\mathcal R}}$. 
To show that ${\bf q} \in Q_A^{\wedge_{\mathcal R}}$ and that 
${\bf p} \equiv {\bf q}$ in $A^{\wedge_{\mathcal R}}$, suppose 
$q = p \cdot a_1 \cdot \ldots \cdot a_n$ for some $a_1, \ldots, a_n \in A$. 
Then \ ${\bf p} \cdot a_1 \cdot \ldots \cdot a_n$ $=$
${\sf red}(r_1 > \ldots > r_{k-1} > p \geq p \cdot a_1 \geq \ \ldots \geq$
$p \cdot a_1 \ldots a_{n-1} \geq p \cdot a_1 \ldots a_{n-1} a_n = q)$ $=$
$(r_1 > \ldots > r_{k-1} > q)$, since $p \equiv q$. Thus, 
${\bf q} \in Q_A^{\wedge_{\mathcal R}}$ and ${\bf p} \equiv {\bf q}$.  

\smallskip

{\sf [Rank condition]} \ Suppose there is a morphism 
$\psi: {\bf A} \to {\bf T}$. Suppose ${\bf q} \cdot \alpha = {\bf q}$ 
in ${\bf A}_A^{\wedge_{\mathcal R}}$, and 
$\eta({\bf q}) \cdot \beta = \eta({\bf q})$ in {\bf A}, and
${\sf rank}_{\bf T}(\alpha) \geq {\sf rank}_{\bf T}(\beta)$. We need to
show that ${\bf q} \cdot \beta = {\bf q}$ in 
${\bf A}_A^{\wedge_{\mathcal R}}$. 
One easily verifies that $\eta({\bf q}) \cdot \beta = \eta({\bf q})$ 
by itself already implies ${\bf q} \cdot \beta = {\bf q}$ (without any
further hypotheses).
 \ \ \ $\Box$

\begin{pro} 
For any automaton {\bf A} in ${\mathcal A}_A$ we have (up to isomorphism):

If \ ${\bf A}_A^{{\sf cov} \ \wedge_{\mathcal R}} = {\bf A}_A^{\sf cov}$ 
 \ then \ ${\bf A}_A^{{\sf cov \, glc} \ \wedge_{\mathcal R}}$ $=$
${\bf A}_A^{\sf cov \, glc}$. 
\end{pro}
{\bf Proof.} Since 
${\bf A}_A^{{\sf cov} \ \wedge_{\mathcal R}} = {\bf A}_A^{\sf cov}$,
the interval \ $[{\bf A}_A^{\wedge_{\mathcal R}},$
  $ {\bf A}_A^{{\sf cov} \ \wedge_{\mathcal R}}]$ \ is contained in
$[{\bf A}, \ {\bf A}_A^{\sf cov}]$, with upper boundaries equal.
The canonical map \, 
${\bf A}_A^{{\sf cov \, glc} \ \wedge_{\mathcal R}}$ 
$\twoheadrightarrow$ ${\bf A}_A^{\sf cov \, glc}$ \, is 
1:1${\mathcal R}$ and satisfies the rank condition (by the previous Prop.).  
Moreover, by definition of $(.)^{\sf glc}$, ${\bf A}_A^{\sf cov \, glc}$
is maximal in the order $\twoheadrightarrow$ for maps that are 
1:1${\mathcal R}$ and that satisfy the rank condition. The Proposition
follows.  
 \ \ \ $\Box$

\bigskip

\noindent {\bf The ${\mathcal L}$-expansion:} \
For any right action $(Q \cdot, S)$ (or any automaton {\bf A} with 
syntactic monoid $S$), we define a {\it left action} of $S$ on 
${\mathcal P}(Q)$ by 

\smallskip

 \ $s \star X = Xs^{-1}$

\smallskip

\noindent for all $s \in S$ and $X \subseteq Q$. Recall the standard notation 
$Xs^{-1} = \{q \in Q: q \cdot s \in X \}$, for any $s \in S$ and 
$X \subseteq Q$.  For a singleton $\{q\}$ we also write $s \star q$ 
and $(q)s^{-1}$ for $s \star \{q\}$, respectively $\{q\}s^{-1}$. 

This is indeed a left action: $(s_1 s_2) \star X$ $=$ $X(s_1 s_2)^{-1}$ 
$=$ $(Xs_2^{-1})s_1^{-1}$ $=$ $s_1 \star (s_2 \star X)$. 
One can prove fairly easily that this is a {\it faithful} left action 
of $S$: If $s_1 \neq s_2$ then, since the given right action is faithful, 
there is $q \in Q$ with $q \cdot s_1 = p_1$ $\neq$ $p_2 = q \cdot s_2$; 
then $q \in (p_1)s_1^{-1}$ but $q \not\in (p_1)s_2^{-1}$.

This action can be restricted to a faithful action on the set
$\{ (q)s^{-1} : \ q \in Q, \ s \in S^1\} \ \subseteq \ {\mathcal P}(Q)$.

The above construction can be applied to automata.
Suppose an automaton ${\bf A} = (Q, i, \cdot)$ has a {\it final} state $f$,
which is reachable from every state in $Q$, i.e., 
$(\forall q \in Q)(\exists w \in A^*)[ \, q \cdot w = f \,]$. Then \   
${\bf A}^{\sf left} = (\{ (f)s^{-1} : s \in S^1\}, \ \{f\}, \ \star)$ 
 \ is a {\it left automaton} with the same syntactic semigroup as {\bf A}. 

\smallskip

For an automaton ${\bf A} = (Q, i, \cdot)$ and its left version 
${\bf A}^{\sf left} =$
$(Q^{\sf left} = \{(f)s^{-1} : s \in S^1\}, \ \{f\}, \ \star)$,
we can construct a left expansion ${\bf A}_A^{\wedge_{\mathcal L}}$;
this generalizes the left Rhodes expansion $S_A^{\wedge_{\mathcal L}}$ 
of a semigroup, and yields a left action (not necessarily faithful) of 
$S_A^{\wedge_{\mathcal L}}$.
The state set of ${\bf A}_A^{\wedge_{\mathcal L}}$ is

\smallskip

$\{ (P_k < \ \ldots \ < P_1 < \{f\}) : k \geq 0, \ {\rm and} $
    $P_k, \ldots, P_1 \in Q^{\sf left} \}$,

\smallskip

\noindent where, for any $P_2, P_1 \in Q^{\sf left}$, we define 

\smallskip

$P_2 \leq P_1$ \ iff \ $(\exists s \in S^1)$ 
  $P_2 = s \star P_1 = P_1 s^{-1}$.  

\smallskip

\noindent The action is defined by 

\smallskip

$a \star (P_k < \ \ldots \ < P_1 < \{f\}) \ = \ $
${\sf red}(a \star P_k \leq P_k < \ \ldots \ < P_1 < \{f\})$.

\smallskip

\noindent The left action of $S_A^{\wedge_{\mathcal L}}$ is 

\smallskip

$(s_m <_{\mathcal L} \ldots <_{\mathcal L} s_1)$ $\star$
$(P_k < \ \ldots \ < P_1 < \{f\}) \ = \ $

\smallskip

${\sf red}(s_m \star P_k \leq \ldots \leq$
  $s_1 \star P_k \leq P_k < \ \ldots \ < P_1 < \{f\})$.
 
\smallskip

\noindent Here, the effect of the reduction operation {\sf red} is to take
the left-most element of an $\equiv$-chain. 

\medskip

When we start with a left action, we can define a corresponding right 
action in a similar way.
This way we can iterate $(.)_A^{\wedge_{\mathcal L}}$ and 
$(.)_A^{\wedge_{\mathcal R}}$. 

\medskip

\begin{pro}
When the automaton expansion $(.)_A^{\wedge_{\mathcal L}}$ is applied to
the right-regular representation $(S^1, 1, \cdot)$ over the alphabet $A$, 
the classical Rhodes expansion $S_A^{\wedge_{\mathcal L}}$ is obtained.
\end{pro}
{\bf Proof.} We saw that $S_A^{\wedge_{\mathcal L}}$ acts on the state 
set of $(S^1, 1, \cdot)_A^{\wedge_{\mathcal L}}$; we want to show that this
action is faithful: If ${\bf s, t} \in S_A^{\wedge_{\mathcal L}}$ act 
in same way on all states of $(S^1, 1, \cdot)_A^{\wedge_{\mathcal L}}$ we 
want to show that ${\bf s} = {\bf t}$.

Let ${\bf s} = $
$(s_m <_{\mathcal L} \ \ldots \ <_{\mathcal L} s_1)$ and ${\bf t} = $
$(t_n <_{\mathcal L} \ \ldots \ <_{\mathcal L} t_1)$.
We will prove by induction that ${\bf s} = {\bf t}$. 
We will use the following (for any $x, y \in S^1$):

\smallskip

\noindent {\sc (l)} \hspace{1.5in} If \ $y <_{\mathcal L} x$ \ then
 \ $(x)y^{-1} = \varnothing$.

\smallskip

\noindent Indeed, if there exists $t \in (x)y^{-1}$ then $ty = x$, which
contradicts $y <_{\mathcal L} x$. Moreover,

\smallskip

\noindent {\sc (i)} \hspace{1.5in} $x = y$ \ iff
 \ $1 \in (x)y^{-1}$.

\smallskip

\noindent For the base-case of the induction we show that $s_n = t_m$ 
follows from \ ${\bf s} \star (\{s_m\})) = {\bf t} \star (\{s_m\})$.
Indeed, ${\bf s} \star (\{s_m\}) = ( (s_m)s_m^{-1} < \ \ldots \ )$ \ $=$ \   
${\bf t} \star (\{s_m\})$ \ $=$  \ $( (s_m)t_n^{-1} < \ \ldots \ )$.
Since $1 \in (s_m)s_m^{-1} = (s_m)t_n^{-1}$, it follows by {\sc (i)} that
$s_m = t_n$.

For the inductive step, let us assume that $s_m = t_n$, $\ \ldots \ $, 
$s_{m-i+1} = t_{n-i+1}$, for some $i>1$; we want to show that 
$s_{m-i} = t_{n-i}$.  We have \

\smallskip
   
${\bf s} \star (\{s_{m-i}\})$ \ $=$ 
 \ ${\sf red}((s_{m-i})s_m^{-1} \leq \ \ldots \ \leq (s_{m-i})s_{m-i+1}^{-1}$
$ \leq (s_{m-i})s_{m-i}^{-1} \leq \ \ldots \ )$ \ $=$

${\bf t} \star (\{s_{m-i}\}) \ = \ $
${\sf red}((s_{m-i})s_m^{-1} \leq \ \ldots \ $
$\leq (s_{m-i})s_{m-i+1}^{-1} \leq (s_{m-i})t_{m-i}^{-1} \leq \ \ldots \ )$. 

\smallskip

\noindent By property {\sc (l)} this becomes 

\smallskip

${\bf s} \star (\{s_{m-i}\}) \ = \ (\varnothing < (s_{m-i})s_{m-i}^{-1} <$
$ \ \ldots \ )$  \ $=$ \ ${\bf t} \star (\{s_{m-i}\}) \ = \ $
$(\varnothing < (s_{m-i})t_{m-i}^{-1} \leq \ \ldots \ )$. 

\smallskip

\noindent Also, $1 \in (s_{m-i})s_{m-i}^{-1} = (s_{m-i})t_{m-i}^{-1}$, hence
$s_{m-i} = t_{m-i}$.

By induction we conclude that \ ${\bf t} =$
$(s_m <_{\mathcal L} \ \ldots \ <_{\mathcal L} s_1 <_{\mathcal L}$
$t_{n-m} <_{\mathcal L} \ \ldots \ <_{\mathcal L} t_1)$, and $m \leq n$.
Moreover, by letting {\bf s} and {\bf t} act on the states 
$(\{t_n\}), \ \ldots \ , (\{t_1\})$, we similarly find that $n \leq m$. 
Hence $n = m$, and ${\bf s} = {\bf t}$.
 \ \ \ $\Box$

\bigskip

\noindent The following can be further developed.

Variations on the definition:
Consider automata with start state and final state.  Redefine the right
and left expansions so that the resulting automata have a start state and 
a final state.

More generally, take automata with a {\sl set} $I$ of initial states and a 
{\sl set} $F$ of final states.  Then ${\bf A}^{\sf left}$ will have 
$Q^{\sf left} = \{ (f)s^{-1} : \, s \in S^1, \, f \in F, \, (\exists t \in
S^1)[ (f)s^{-1}t^{-1} \cap I \neq \varnothing ]\}$, and it will also have 
a set of initial states $I^{\sf left} = \{ \{f\} : f \in F\}$ and a set of
final states $F^{\sf left} =$
$\{P \in Q^{\sf left} : P \cap I \neq \varnothing\}$.  And both 
${\bf A}_A^{\wedge_{\mathcal L}}$ and ${\bf A}_A^{\wedge_{\mathcal R}}$ 
will have a set of initial states and a set of final states.

Now we can iterate the left and right expansions and obtain two-sided 
unambiguity.


\section{A bottom-up inductive construction of the glc cover}

\subsection{Rank functions}

The rank functions that we defined earlier will now be extended to edges 
and to walks of an $A$-automaton, within a reachability class. 
Let ${\bf A} = (Q, \cdot)$ be a finite automaton, let 
$\psi: {\bf A} \to {\bf B}$ be an $A$-automaton morphism. 
An {\em edge} is of the form $(q, \, a, \, q \cdot a)$, where 
$q \in Q, a \in A$; we will also simply write $(q,a)$. We say that an edge 
is in a reachability class iff $q \equiv_{\mathcal R} q \cdot a$.

More  generally, a {\em walk} of length $n$ is of the form
$(q, \, a_1, \, q_1, \, a_2, \, q_2, \ \ldots \ , a_i, \, q_i, \, a_{i+1},$ 
$ \ \ldots \ , a_n, \, q_n)$, where $q, q_1, \ldots, q_n \in Q$ and 
$a_1, a_2, \ldots , a_n \in A$, with  $q_i = q a_1 \ldots a_i$ for 
$i = 1, 2, \ldots n$. We will also simply denote a walk by 
$(q, \, a_1 \ldots a_n)$.
A walk $(q, \, a_1 \ldots a_n)$ is said to be {\em within a reachability 
class} \, iff \, $q \cdot a_1 \ldots a_n \equiv_{\mathcal R} q$.

\begin{defn} \label{edgerank} \   
When $q \equiv_{\mathcal R} q \cdot a$ we define the {\em edge-rank} of
$(q, a)$ by \

\smallskip

${\sf rank}_{\bf B}(q,a)$ $ \ = \ $
${\sf min} \{ {\sf rank}_{\bf B}(a x) : \ x \in A^*, \ q \cdot a x = q\}$. 

\smallskip
 
\noindent For a path $(q,w)$ within a reachability class we define the
{\em path-rank} by 

\smallskip

${\sf rank}_{\bf B}(q,w)$ $ \ = \ $
${\sf min} \{ {\sf rank}_{\bf B}(w x) : \ x \in A^*, \ q \cdot w x = q\}$.
\end{defn}


\noindent {\bf Notation:} \ For an $A$-automaton ${\bf A} = (Q, \cdot)$ and 
a word $w \in A^*$, $[w]$ denotes the element represented by $w$ in the 
syntactic monoid of {\bf A}. This means that $[w]$ is the congruence class 
of $w$ for the monoid congruence $\equiv_{\bf A}$ on $A^*$ defined by 
 \, $x \equiv_{\bf A} y$ \, iff \, 
$(\forall q \in Q)[ \, q \cdot x = q \cdot y \, ]$.


\begin{defn} 
 \ An $A$-automaton {\bf A} has {\em idempotent stabilizers} iff for 
every state $q$ and every word $w \in A^*:$ \ $q \cdot w = q$ implies 
$[w]^2 = [w]$.
\end{defn}

\begin{defn}
 \ An $A$-automaton {\bf A} has {\em ${\cal R}$-trivial stabilizers} iff 
for every state $q$ and every word $w \in A^*:$ \ $q \cdot w = q$ implies
that the $\equiv_{\mathcal R}$-class of $[w]$ is a singleton (in the 
syntactic monoid of {\bf A}).  
\end{defn}

\begin{pro}  \label{rankInvarConjug} 
{\em (Invariance of rank under conjugation)} \  Let ${\bf A} = (Q, \cdot)$
be an $A$-automaton with idempotent stabilizers, and suppose $(q, w)$ 
is a closed walk, i.e., $q \cdot w = q$, $w \in A^*$, $q \in Q$. 
If $w$ is factored as $w = xy$, then 
 \ ${\sf rank}_{\bf B}(xy) = {\sf rank}_{\bf B}(yx)$. 
\end{pro}
{\bf Proof.} Since $q w = q$ and stabilizers are idempotents, we have 
$[w]^2 = [w]$. Also, $q w = q$ implies $q w x = q x$, i.e., 
$q x y x = q x$ (since $w = xy$). So, $yx$ stabilizes $qx$, hence 
$[yx]^2 = [yx]$. 
Moreover, $[xy] \equiv_{\mathcal J} [yx]$; indeed, 
$[xy] = [xyxy] \leq_{\mathcal J} [yx]$, and 
$[yx] = [yxyx] \leq_{\mathcal J} [xy]$.
Now, since $[xy]$ and $[yx]$ are $\equiv_{\mathcal J}$-related 
idempotents, they have the same algebraic rank.   \ \ \ $\Box$

\bigskip

For Prop.\ \ref{pathRank} we will need an expansion with special properties;
we first make sure that such an expansion is available. Let 
$(.)_A^{\wedge_{\cal L}}$ be the Rhodes expansion, let 
$(.)^{\sf RB}$ be the rectangular-bands expansion {\bf RB} \mal $(.)$,
and let $(.)^{\sf IS}$ be the expansion from \cite{IdempotStab}.
For properties of the Rhodes expansion, see \cite{RhInfItII} Appendix A.IV,
\cite{ItExp}, and Tilson's chapter XII in \cite{EilenbergB}.
We will use the fact that \  
$(.)^{{\sf RB} \, \wedge_{\cal L}} = (.)^{\sf RB}$, i.e., 
$(.)^{\sf RB}$ is stable under $(.)^{ \wedge_{\cal L}}$. 
 
\begin{pro}  
For any $A$-generated semigroup $S_A$, the expansions 
$S_A^{ {\sf IS} \, \wedge_{\cal L} }$ and $S_A^{\sf IS \, RB}$ in the 
category of $A$-generated semigroups have the following properties:
All right-stabilizers are $\cal R$-trivial bands, and the $\cal L$-order
in each right-stabilizer is unambiguous.
\end{pro}
{\bf Proof.} In \cite{IdempotStab} it is proved that in $S_A^{\sf IS}$ the 
right-stabilizers are $\cal R$-trivial bands. For any $A$-generated 
semigroup $T_A$ the natural maps 
$\eta_{\cal L}: T_A^{\wedge_{\cal L}} \to T_A$ and
$\eta_{\sf RB}: T_A^{\sf RB} \to T_A$ have the property that the inverse 
image $\eta^{-1}(B)$ of any band $B \subseteq T_A$ is a band (where 
$\eta$ stands for either $\eta_{\cal L}$ or $\eta_{\sf RB}$).  
And any stabilizer $\Sigma' \subseteq T_A^{\wedge_{\cal L}}$ is contained
in $\eta^{-1} \eta(\Sigma')$, and $\eta(\Sigma')$ is contained in a stabilizer
$\Sigma$ in $T_A$; so any stabilizer $\Sigma' \subseteq T_A^{\wedge_{\cal L}}$
is contained in $\eta^{-1}(\Sigma)$ where $\Sigma$ is a stabilizer of $T_A$ . 
Hence any stabilizer in $S_A^{ {\sf IS} \, \wedge_{\cal L} }$ or in
$S_A^{\sf IS \,  RB}$ is a band.

Moreover, for any $T_A$, every right-stabilizer in $T_A^{\wedge_{\cal L}}$
is $\cal R$-trivial and the $\cal L$-order of $T_A^{\wedge_{\cal L}}$ is 
unambiguous. 
Since $T_A^{\sf RB} = T_A^{{\sf RB} \, \wedge_{\cal L}}$, 
every right-stabilizer is $\cal R$-trivial there too, and the $\cal L$-order 
is unambiguous.  
 \ \ \ $\Box$

\begin{pro}  \label{pathRank}           
{\bf (Path-rank vs.\ edge-rank).} 
 \ Let {\bf A} be an $A$-automaton whose stabilizers are $\cal R$-trivial
bands with unambiguous $\cal L$-order. 
Let $(q_0, \, a_1 \ldots a_n)$ be a walk within an ${\mathcal R}$-class, 
i.e., $q_0 a_1 \ldots a_n \equiv_{_{\mathcal R}} q_0$. Let $q_i$ denote 
$q_0 a_1 \ldots a_i$ for $i = 1, \ldots, n$. Then the rank of the path 
is equal to the maximum of the edge-ranks, i.e., 

\smallskip

\noindent {\small \rm (1)} \hspace{.6in}
${\sf rank}_{\bf B}(q_0, \, a_1 \ldots a_n) \ = \ $
${\sf max} \{ \, {\sf rank}_{\bf B}(q_{i-1}, a_i) : i = 1, \ldots, n\}$.

\smallskip

\noindent If the path $(q_0, \, a_1 \ldots a_n)$ is closed (i.e.,
$q_0 a_1 \ldots a_n = q_0$) then the rank of the path is equal to the 
rank of the labeling word, i.e., 

\smallskip

\noindent {\small \rm (2)} \hspace{.6in}
  \ ${\sf rank}_{\bf B}(q_0, \, a_1 \ldots a_n) \ = \ $
 ${\sf rank}_{\bf B}(a_1 \ldots a_n)$.  
\end{pro}
(This Proposition will also be called the ``{\bf Sausage Lemma}'' since a 
path between several reachability classes can be pictured as a string of 
sausages.) 

\noindent
{\bf Proof.} In the proof we drop the ubiquitous subscript {\bf B} of 
{\sf rank}. 

\smallskip

\noindent [(1) $\geq$] \ By the definition of path-rank there exists 
$\beta \in A^*$ such that ${\sf rank}(q_0, \, a_1 \ldots a_n)$ $ = $
${\sf rank}(a_1 \ldots a_n \beta)$ 
and $q_0 a_1 \ldots a_n \beta = q_0$.
The word $\beta$ is the ``$x$'' in the definition of path-rank for which the 
{\sf min} is achieved.  Then we have for any $j = 1, \ldots , n$:  

\smallskip

${\sf rank}(q_0, \, a_1 \ldots a_n) \ = \ $
${\sf rank}(a_1 \ldots a_n \, \beta) \ $

 \ \ \ $= \ {\sf rank}(a_j \, a_{j+1} \ldots a_n \beta a_1 \ldots a_{j-1})$  
 \hspace{.5in} (by Prop.\ \ref{rankInvarConjug})  

\ \ \ $\geq \ {\sf min} \{ {\sf rank}(a_j x) : $
    $x \in A^*, \ q_{j-1} a_j x = q_{j-1}\}$   \hspace{.5in} (since 
    $q_{j-1} \cdot a_j \ldots a_n \beta a_1 \ldots a_{j-1} $
    $= q_{j-1}$) 

 \ \ \ $= \ {\sf rank}(q_{j-1}, a_j)$ \hspace{.5in} 
                                  (by the definition of edge-rank).

\smallskip

\noindent Therefore, 
${\sf rank}(q_0, a_1 \ldots a_n) \geq {\sf rank}(q_{j-1}, a_j)$ 
for every $j = 1, \ldots , n$; hence,
${\sf rank}(q_0, a_1 \ldots a_n) \geq {\sf max}\{ {\sf rank}(q_{j-1}, a_j):$
$j = 1, \ldots , n\}$. 

\smallskip

\noindent [(1) $\leq$] \ We use induction on the path-length $n$.
For $n = 1$, the two sides of (1) are identical.
For $n > 1$, assume \, ${\sf rank}(q_0, a_1 \ldots a_{n-1}) \ \leq \ $
$ {\sf max}\{{\sf rank}(q_{j-1}, a_j): j = 1, \ldots, n-1 \}$;
equivalently, ${\sf rank}(q_0, a_1 \ldots a_{n-1}) \ \leq \ $
${\sf rank}(q_{j-1}, a_j)$ for some $j = 1, \ldots, n-1$.
We want to prove that ${\sf rank}(q_0, a_1 \ldots a_{n-1} a_n) \leq $
${\sf rank}(q_{i-1}, a_i)$ for some $i = 1, \ldots, n-1, n$. 

Let $\beta \in A^*$ be such that 
 \, $q_0 a_1 \ldots a_{n-1} \beta = q_0$ \, and \, 
${\sf rank}(q_0, a_1 \ldots a_{n-1}) = {\sf rank}(a_1 \ldots a_{n-1} \beta)$,  
i.e., the {\sf min} in the definition of path-rank is reached when $x$ is 
$\beta$.  Similarly, let $\beta' \in A^*$ be such that 
 \, $q_{n-1} a_n \beta' = q_n$ \, and \, 
${\sf rank}(q_{n-1}, a_n) = {\sf rank}(a_n \beta')$. 
Then each of the following is a closed path: 

$(q_{n-1}, \ a_n \beta')$, \ \
$(q_{n-1}, \ \beta a_1 \ldots a_{n-1})$, \ \    
$(q_0, \ a_1 \ldots a_{n-1} \beta)$ , \ and \
$(q_0, \  a_1 \ldots a_{n-1} a_n \beta' \beta)$.

\smallskip

\noindent It follows that 
 \ $[a_n \beta']$,  
 \ $[\beta a_1 \ldots a_{n-1}]$,
 \ $[a_1 \ldots a_{n-1} \beta]$, \ and
 \ $[a_1 \ldots a_{n-1} a_n \beta' \beta]$
 \ belong to stabilizers, hence they are idempotents. 

By the definition of rank (using min), 
${\sf rank}(q_0, \, a_1 \ldots a_{n-1} a_n) \leq $
${\sf rank}(a_1 \ldots a_{n-1} a_n \beta' \beta)$, and by Prop.\ \ref{rankInvarConjug},
 \, ${\sf rank}(a_1 \ldots a_{n-1} a_n \beta' \beta) \ = \ $
${\sf rank}(\beta a_1 \ldots a_{n-1} a_n \beta')$.

\smallskip

\noindent {\sf Case 1:} \ $[\beta a_1 \ldots a_{n-1}]$ $\equiv_{\cal R}$
 $[\beta a_1 \ldots a_{n-1} a_n \beta']$

Then by $\cal R$-triviality of stabilizers, 
$[\beta a_1 \ldots a_{n-1} a_n \beta'] = [\beta a_1 \ldots a_{n-1}]$, and 
these are idempotents. Hence,
${\sf rank}(q_0, \, a_1 \ldots a_{n-1} a_n) \leq $
${\sf rank}([\beta a_1 \ldots a_{n-1} a_n \beta']) = $
${\sf rank}([\beta a_1 \ldots a_{n-1}]) = $ 
${\sf rank}(q_0, \, a_1 \ldots a_{n-1}) \leq {\sf rank}(q_{j-1}, a_j)$
for some $j$ (the latter by the induction hypothesis). 

\smallskip

\noindent {\sf Case 2:} \ $[\beta a_1 \ldots a_{n-1}]$ $>_{\cal R}$ 
$[\beta a_1 \ldots a_{n-1} a_n \beta']$
$\leq_{\cal L}$ $[a_n \beta']$.

\smallskip

\noindent {\sf Case 2.1:} \  
$[\beta a_1 \ldots a_{n-1} a_n \beta'] \equiv_{\cal L} [a_n \beta']$

Then, since these are idempotents,
${\sf rank}(q_0, \, a_1 \ldots a_{n-1} a_n) \leq $
${\sf rank}([\beta a_1 \ldots a_{n-1} a_n \beta']) = $
${\sf rank}([a_n \beta']) = {\sf rank}(q_{n-1}, a_n)$.

\smallskip

\noindent {\sf Case 2.2:} \ $[\beta a_1 \ldots a_{n-1}]$ $>_{\cal R}$
$[\beta a_1 \ldots a_{n-1} a_n \beta'] <_{\cal L} [a_n \beta']$.

Let $e_1 = [\beta a_1 \ldots a_{n-1}]$, \ $e_2 = [a_n \beta']$; 
then $e_1 e_2 = [\beta a_1 \ldots a_{n-1} a_n \beta']$, and 
$e_1 e_2 = (e_1 e_2)^2$, and $e_1 >_{\cal R} e_1 e_2 <_{\cal L} e_2$. 
Since stabilizers are $\cal R$-trivial bands, $e_1 e_2 = e_1 e_2 e_1$. 
Similarly, $e_2 e_1 = e_2 e_1 e_2$. 

Thus, $e_1 >_{\cal L} e_2 e_1 = e_2 e_1 e_2 <_{\cal L} e_2$, which by 
unambiguity of the $\cal L$-order implies that either $e_1 \leq_{\cal L} e_2$,
i.e., $e_1 = e_1 e_2$, or $e_2 \leq_{\cal L} e_1$, i.e., $e_2 = e_2 e_1$. 
But $e_1 >_{\cal R} e_1 e_2$ contradicts $e_1 = e_1 e_2$.
And $e_1 e_2 <_{\cal L} e_2$ contradicts $e_2 = e_2 e_1$ (since 
$e_1 e_2 = e_1 e_2 e_1 \equiv_{\cal L} e_2 e_1$).
Thus, when the syntactic monoid of the automaton {\bf A} has unambiguous 
$\cal L$-order case 2.2 is impossible. 

\medskip

\noindent Proof of (2): 
 \ By definition, ${\sf rank}(q_0, a_1 \ldots a_n) =$
${\sf min}\{ {\sf rank}(a_1 \ldots a_n x) : q_0 a_1 \ldots a_n x = q_0\}$.
We assumed that $a_1 \ldots a_n$ stabilizes $q_0$. 
By the properties of stabilizers, $[a_1 \ldots a_n]$ and $[a_1 \ldots a_n x]$ 
are idempotents.  Since $a_1 \ldots a_n x \leq_{\cal J} a_1 \ldots a_n$ 
and since these are idempotents, we have 
${\sf rank}(a_1 \ldots a_n x) \geq {\sf rank}(a_1 \ldots a_n)$ 
for any $[x]$, with equality when $x$ is empty. Hence, 
${\sf min}\{ {\sf rank}(a_1 \ldots a_n x) : q_0 a_1 \ldots a_n x = q_0\}$
$= {\sf rank}(a_1 \ldots a_n)$. 
 \ \ \ $\Box$

\bigskip

For the rest of this subsection let $T_A$ be an $A$-generated semigroup.  
Every $\cal R$-class $R$ of $T_A$ will be viewed as a deterministic partial
finite automaton (with state set $R$ and alphabet $A$). This is a strongly 
connected component of the right Cayley graph $\Gamma(T_A)$. 
We call this the $\cal R$-class automaton of $R$.

We assume that the semigroup $T_A$ satisfies the assumptions of Prop.\ 
\ref{pathRank}, i.e., that its stabilizers are $\cal R$-trivial bands 
with unambiguous $\cal L$-order.

We consider the algebraic rank of elements of $T_A$ with respect to the
identity map on $T_A$, denoted by ${\sf rank}(.)$.  
We also consider the edge-rank function, as in Definition \ref{edgerank}, 
with respect to the identity map on $T_A$. We also call this rank function 
${\sf rank}(.)$. So for an edge $(q, a)$ with $q, qa \in R$ we have:
${\sf rank}(q, a) = {\sf min}\{ {\sf rank}(ax) : x \in A^*, \ qax = a\}$;
here, ${\sf rank}(ax)$ is the algebraic rank.

\begin{defn} 
For the $\cal R$-class $R$ of an $A$-automaton we define 

\smallskip

 \ \ \ ${\sf rank}(R) \ = \ $
${\sf max}\{ {\sf rank}(q, a) \, : \, q, \, qa \in R, \, a \in A \}$.  

\smallskip

\noindent For $r \in R$ and $j = 0, 1, \ldots$ we define 

\smallskip

 \ \ \ ${\sf paths}(r, j) \ = \ $
$\{ w \in A^* : \, r \cdot w \in R$ and every edge $e$ of the 
path $(r, w)$ in the automaton $R$ 

\hspace{1.8in} has rank ${\sf rank}(e) \leq j\}$,

\smallskip

 \ \ \ $P(r, j)$ ($\subseteq R$) is the set of vertices of the paths in 
${\sf paths}(r, j)$,

\smallskip

 \ \ \ $E(r, j)$ is the set of edges of the paths in ${\sf paths}(r, j)$.
\end{defn}

\begin{pro} 
For every $\cal R$-class $R$ of an $A$-automaton and every fixed 
$j_o = 0, 1, \ldots$  we have: The set of sets 
$\{ P(r, j_o) : r \in R\}$ is a partition of $R$.
\end{pro}
{\bf Proof.} Since $r \in P(r, j_o)$ the sets $P(r, j_o)$ are obviously 
non-empty and $R = \, \bigcup_{r \in R} P(r, j_o)$.
Let us prove that non-equal sets that overlap are equal. If 
$r \in P(r_i, j_o)$, then there is a path in 
$R$ from $r_i$ to $r$, labelled by some $\alpha \in A^+$, with maximum
edge rank $\le j_o$. Then by Prop.\ \ref{pathRank} (Sausage Lemma), this 
path can be extended to a cycle containing $r_i$ with no increase in 
maximum edge rank.  Hence, if $r \in P(r_1, j_o) \cup P(r_2, j_o)$ then 
there are paths from $r_1$ to $r_2$ and from $r_2$ to $r_1$ of maximum
edge rank $\le j_o$; so $P(r_1, j_o) = P(r_2, j_o)$. 
 \ \ \ $\Box$

\bigskip

\noindent {\bf Examples of $P(.,.)$ and ${\sf paths}(.,.)$:}  

\smallskip

\noindent {\bf (1)} \ Suppose the state diagram of $R$, with ranks above 
and below the edges, is 

%

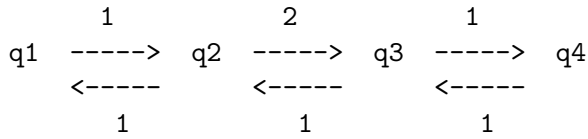
\begin{verbatim}
          1           2           1
    q1  ----->  q2  ----->  q3  ----->  q4
        <-----      <-----      <-----
           1           1           1
\end{verbatim}

\noindent Then the partition for $j = 0$ is 
$\, \{ \{q_1\}, \{q_2\}, \{q_3\}, \{q_4\} \}$, for $j = 1$ it is 
$\, \{ \{q_1, q_2\}, \{q_3, q_4\} \}$, and for $j = 2$ it is
$\, \{ \{q_1, q_2, q_3, q_4\} \}$. 

\medskip

\noindent {\bf (2)} \ Suppose the state diagram of $R$ consists of 
the states $q_1$ and $q_2$, with edges $(q_1, a \, (1), q_2)$, \,
$(q_1, b \, (2), q_2)$, and $(q_2, c \, (1), q_1)$, with ranks in 
parentheses.
Then ${\sf paths}(q_1, 1)$ consists of all paths from $q_1$ with label in
$(ac)^* \cup (ac)^*a$. 
And $P(q_1, 1) = P(q_2, 1) = P(q_1, 2) = P(q_2, 2) = \{ \{ q_1, q_2 \} \}$. 

\noindent {\sf [End, Examples.]}

\bigskip

We are now going back to the covers
 \ $U_A = T_A^{\sf cov \ IS \ RB} \ $
$\stackrel{\varphi_1}{\twoheadrightarrow} \ [T_A, U_A]^{\sf glc}$
$ \ \stackrel{\varphi_2}{\twoheadrightarrow} \ T_A$, where $\varphi_2$ is
1:1$\cal R$ and satisfies the rank condition.  
In this setting we will also write {\sf UP} for $U_A$, {\sf MIDDLE} for 
$[T_A, U_A]^{\sf glc}$, and {\sf DOWN} for $T_A$.
Let $\tilde{R}$ be an $\cal R$-class of {\sf MIDDLE}, let $R$ be the 
$\cal R$-class of {\sf DOWN} such that 
$\varphi_2: \tilde{R} \ \hookrightarrow \ R$. 
To assign a rank to an edge 
$\tilde{r} \stackrel{a}{\longrightarrow} \tilde{r} a$ of $\tilde{R}$, we 
apply $\varphi_2$, yielding $(\tilde{r})\varphi_2$
$\stackrel{a}{\longrightarrow} (\tilde{r} a)\varphi_2$, and take the rank of
this arrow in $R$. Recall that ${\sf rank}(\tilde{R})$ is the max rank over
all idempotents in $E(\tilde{R})$ ($= E(R)$). Notice that by Prop.\
\ref{pathRank}, ${\sf rank}(R)$ is the same as algebraic rank of all 
idempotents that fix any element of $R$ under right-multiplication.
We also have 
$\tilde{n} = {\sf rank}(\tilde{R}) \le {\sf rank}(R) = n$.

For an $\cal R$-class $C$ of a semigroup, viewed as an $A$-automaton, 
we denote the edge set of $C$ by $E(C)$. 
As above, $E(r, j)$ is the set of edges of ${\sf paths}(r,j)$ in the 
$\cal R$-class automaton $R$. 
 
\begin{pro} \label{edgesPaths} \   
For $\tilde{R}$ and $R$ as above and for any $\tilde{r} \in \tilde{R}$ and
$r = (\tilde{r})\varphi_2 \in R$ we have:

\smallskip
 
 \ \ \ \ \ $E(\tilde{R}) \ = \ E(r, \tilde{n}) \ \subseteq \ E(R)$.  
\end{pro}
{\bf Proof.}  Let $\alpha \in A^*$ be the label of a cycle in $\tilde{R}$,
starting at $\tilde{r}$, that traverses every edge of $\tilde{R}$. 
Then each path in ${\sf paths}(r, \tilde{n})$ can be extended to a cycle 
that returns to $r$, and that uses only edges of edge-rank $\le \tilde{n}$
(by Prop.\ \ref{pathRank}, the ``Sausage Lemma'').
Now the rank condition implies the result as follows: 
We have $\tilde{r} \alpha = \tilde{r}$, and $\alpha$ has at least one edge 
of rank $\tilde{n}$, and passes through all vertices of $\tilde{R}$. 
If $r \beta = r$ and ${\sf rank}(\beta) \le {\sf rank}(\alpha) = \tilde{n}$
then $\tilde{r} \beta = \tilde{r}$, by the rank condition. So, 
$E(\tilde{R}) = E(r, \tilde{n}) \subseteq E(R)$.  
 \ \ \ $\Box$ 

\begin{cor}
If ${\sf rank}(\tilde{R}) = {\sf rank}(R)$ then the restriction of $\varphi_2$ 
to $\tilde{R}$ is an isomorphism from $\tilde{R}$ to $R$.  
 \ \ \  $\Box$
\end{cor}
We will give a formula for ${\sf rank}(R)$ in terms of idempotents.
The following Lemma will be used. 
For idempotents $e, f$ we write \, $e \ge f$  \ iff \ $f = ef = fe$; 
this is called the {\em idempotent order}, or equivalently,
the $\cal H$-order between idempotents; it is a partial order. We write
$e > f$ iff $e \ge f$ and $e \ne f$.   

\begin{lem} \label{regJclassIdemp} {\rm (Rhodes 1966).} \ Let $J_1$ and 
$J_2$ be regular $\cal J$-classes of a finite semigroup, such that 
$J_1 >_{\cal J} J_2$. Then for every idempotent $e \in J_1$ there 
exists an idempotent $f \in J_2$ such that $e > f$.
\end{lem}
{\bf Proof.} See the proof of Prop.\ 3.1 in \cite{Rh66}. 
  \ \ \  $\Box$

\medskip

\begin{pro} \label{rankRwithIdemp} \
Suppose $R$ is a regular $\cal R$-class of $T_A$ ($= {\sf DOWN}$) such that
the right-stabilizer of any element of $R$ is an $\cal R$-trivial band 
with unambiguous $\cal L$-order.
Then ${\sf rank}(R) = n$ implies that there exist $r \in R$, and 
idempotents $e_0, e_1, \ \ldots \ , e_n$, such that 

\smallskip
 
 \ \ \  \ \ \   
 $e_0  > e_1 > \ \ldots \ > e_n \equiv_{\cal L} r$ ,

\smallskip
 
\noindent where ``$>$'' is the strict idempotent order.  Hence the stabilizer
$(r){\sf St}$ contains the idempotent-order chain 
$e_0  > e_1 > \ \ldots \ > e_n$. 

\end{pro}
{\bf Proof.} By definition, ${\sf rank}(R)$ is the maximum of the algebraic
ranks of all idempotents that fix any element of $R$ under right
multiplication. So if ${\sf rank}(R) = n$, there exists $e = e^2 \in R$ with 
${\sf rank}(e) = n$, such that 
 \ $e \equiv_{\cal J} e'_n <_{\cal J} \ \ldots \ <_{\cal J} e'_1 <_{\cal J}$
$e'_0$ \, (a $\cal J$-chain of idempotents). 
Applying Lemma \ref{rankRwithIdemp} $n$ times yields $e \equiv_{\cal J} $
$e_n < \ \ldots \ < e_1 < e_0$, a chain in the idempotent order.
Hence there exists $r \in R$ with $e \equiv_{\cal R} r$ and 
$r \equiv_{\cal L} e_n$. So, $(r){\sf St} = (e_n){\sf St}$, which contains 
$e_n < \ \ldots \ < e_1 < e_0$.  
  \ \ \   \ \ \  $\Box$

\bigskip
 
\noindent {\bf Examples.} \ \ \ \ \ \ \

\smallskip

\noindent The following examples show that the stabilizers of the elements
$r \in R$ are not all isomorphic. And when ${\sf rank}(R) = n$,
some of these stabilizers contains an idempotent chain 
$e_n < \ \ldots \ < e_1 < e_0$, and some stabilizers do not.  

\smallskip

\noindent {\bf (a)} \ Consider the semigroup $T = \{f, a, b\}$ with the
multiplication table

\bigskip

\hspace{1in}
\begin{tabular}{r||r|r|r}
$\bullet$      & $f$ & $a$ & $b$  \\ \hline \hline
$f$            & $f$ & $a$ & $b$  \\ \hline
$a$            & $b$ & $a$ & $b$  \\ \hline
$b$            & $b$ & $a$ & $b$  \\   
\end{tabular} 

\bigskip

\noindent Then the rank of the ${\cal R}$-class $\{a, b\}$ is 1,
$(a){\sf St} = \{a\}$, and $(b){\sf St} = \{b, f\}$ (with $b < f$). 
So $(b){\sf St}$ contains a chain of idempotents of length 1 (counting
the inequalities in the chain), while $(a){\sf St}$ does not.  

\medskip

\noindent {\bf (b)} \ Suppose the following diagram is part of the Cayley
graph of a semigroup $S_A$. The vertices are $\{q_1, q_2, q_3, q_4\}$ and the
edges are

$(q_1, a, q_2)$, \, $(q_2, b, q_1)$, \, both with rank 1;

$(q_2, c, q_3)$, \, $(q_3, d, q_2)$, \, both with rank 2;

$(q_3, e, q_4)$, \, $(q_4, f, q_3)$, \, both with rank 3.

\smallskip

\noindent Then under the hypotheses of Prop.\ \ref{pathRank}, the ranks 
of the right-stabilizers are as follows.

$(q_1){\sf St} \,$ has elements of ranks 1, 2, 3;

$(q_4){\sf St} \,$ has rank 3 only;

$(q_2){\sf St} \,$ has ranks 1, 2, 3.;

$(q_3){\sf St} \,$ has ranks 2 and 3.

\smallskip

\noindent This is proved as follows. 

First, if a word $\ell_j \in \{a,b,c,d,e,f\}^+$ (for $j = 1, 2$) 
fixes $q \in \{q_1, q_2, q_3, q_4\}$, i.e., $q\ell_j = q$, then $\ell_j$
is an idempotent (by the hypotheses of Prop.\ \ref{pathRank}); moreover, 
$\ell_2 \ell_1 = \ell_2$ iff $\ell_2 \le_{\cal L} \ell_1$.
Also, the $\cal L$-order of $(q){\sf St}$ is the $\cal L$-order of the whole
semigroup since $\ell_2, \ell_1$ are idempotents. Since $(q){\sf St}$ has
unambiguous $\cal L$-order, either $\ell_2 \le_{\cal L} \ell_1$ or
$\ell_1 \le_{\cal L} \ell_2$; i.e., either $\ell_2 \ell_1 = \ell_2$ or
$\ell_1 \ell_2 = \ell_1$. The rank condition determines the direction of 
the $\cal L$-order: ${\sf rank}(\ell_2) \ge {\sf rank}(\ell_2)$ implies
$\ell_2 \ell_1 = \ell_2$, etc., as is easy to see. 

Second, we have already proved the conjugation property of loops (Prop.\
\ref{rankInvarConjug}), i.e., the rank of a loop does not depend on the
chosen starting point on the loop. We will denote conjugation by $\sim$. 

 From these two observations we now prove  the rank properties of the
stabilizers.
To see that ${\sf St}(q_1)$ has ranks 1, 2, 3, we consider the loops
$(q_1, ab)$, \ $(q_1, acdb)$, \ $(q_1, acefdb)$ with start point $q_1$; 
these loops have rank respectively 1, 2, 3.
For example, to see that $acefdb$ has rank 3, we observe that
$(q_1, acefdb) \sim (q_3, efdbac)$, and 
$(q_3, dbac) \sim (q_2, cdba)$; moreover, $cd$ has rank 2, and $ba$ has 
rank 1. 
So, ${\sf rank}(cdba) = {\sf rank}(cd)$, and 
$dbac \sim cdba$, $cd \sim dc$; so ${\sf rank}(dbac) = {\sf rank}(dc) = 2$,
Since ${\sf rank}(ef) = 3$ and ${\sf rank}(dbac) = 2$, we have
${\sf rank}(efdbac) = {\sf rank}((ef)(dbac)) = {\sf rank}(ef)  = 3$.
So, ${\sf rank}(acefdb) = 3$.

\bigskip

%
%
%
%
%
%
%
%
%
%
% ... 
%
%
%
%
% ...
%
%

\bigskip

\subsection{Bottom-up construction}

We construct the cover ${\bf A}^{\sf cov \, glc}$ in a different way 
than before, using ``bottom-up'' induction. We will later prove that under 
certain conditions this {\it bottom-up construction} yields the same 
cover ${\bf A}^{\sf cov \, glc}$ as defined earlier. However, in the 
meanwhile we need to distinguish the two, and we will denote the cover
resulting from the bottom-up construction by ${\bf A}^{\sf cov \, GLC}$. 

We start with a cover morphism 
$\varphi: {\bf A}^{\sf cov} \to {\bf A}$ and the interval 
$[{\bf A}, {\bf A}^{\sf cov}]$, which is a finite lattice.  
Since ${\bf A}^{\sf cov \, GLC}$ is intended to belong to the 
interval $[{\bf A}, {\bf A}^{\sf cov}]$, we want to construct morphisms 
$\varphi_1: {\bf A}^{\sf cov} \to {\bf A}^{\sf cov \, GLC}$ and
$\varphi_2: {\bf A}^{\sf cov \, GLC} \to {\bf A}$ such that 
$(.)\varphi = (.)\varphi_1 \circ \varphi_2$, and such that $\varphi_2$ is
1:1${\mathcal R}$ and satisfies the rank condition.
Moreover, $\varphi_2$ should be maximal (i.e., ${\sf mod} \varphi_1$ should 
be maximally fine) with respect to these properties.
Hence the congruence ${\sf mod} \varphi_1$ is a refinement of the congruence 
${\sf mod} \varphi$ on the state set $Q^{\sf cov}$ of ${\bf A}^{\sf cov}$, 
and every ${\sf mod} \varphi_1$ class is mapped to one element of $Q$ by 
$\varphi$ (where $Q$ is the state-set of {\bf A}). We will write states of 
${\bf A}^{\sf cov}$ with an overline to make them recognizable.  
For a state ${\overline q}$ of ${\bf A}^{\sf cov}$, we denote the 
${\sf mod} \varphi_1$ congruence class of ${\overline q}$ by 
$[{\overline q}]_{{\sf mod} \varphi_1}$ or more briefly by
$[{\overline q}]_{\varphi_1}$.

We will construct ${\bf A}^{\sf cov \, GLC}$ by defining the congruence 
${\sf mod} \varphi_1$ on $Q^{\sf cov}$, and by defining $\varphi_2$ on each 
constructed congruence class. 
The construction of $[{\overline q}]_{\varphi_1}$ proceeds by {\it induction 
on the directed path-length from ${\overline \i}$ to} ${\overline q}$ in 
${\bf A}^{\sf cov}$, where ${\overline \i}$ is the start state of 
${\bf A}^{\sf cov}$.  
In this induction we assume that the $\cal R$-order of each 
${\bf A}^{\sf cov}$ is a tree; this will hold if $(.)^{\sf cov}$ is closed 
under the right Rhodes expansion $(.)^{\wedge_{\cal R}}$. 

\medskip

First, the congruence class $[{\overline \i}]_{\varphi_1}$ is 
$\{{\overline \i}\}$, since the ${\sf mod} \varphi$ congruence class of 
${\overline \i}$ is just $\{{\ov \i}\}$ (assuming that ${\ov \i}$ 
is not reachable from any state).

\smallskip

For the {\it inductive step} we assume that for every state $\ov{q}_1$ with
a certain directed path-length from ${\overline \i}$, a congruence class 
$[{\overline q}_1]_{\varphi_1}$
$ = \{ {\overline q}_1, \ \ldots \ , {\overline q}_{\ell}\}$ 
has been constructed. Our goal is to construct the congruence class 
$[{\overline q}_1 \cdot a]_{\varphi_1}$ for each $a \in A$ such that
the directed path-length from ${\overline \i}$ to $\ov{q}_1 \cdot a$
is larger than the directed path-length from ${\overline \i}$ to
$\ov{q}_1$. 

Note that $[{\overline q}_j \cdot a]_{\varphi_1}$ $ = $
$[{\overline q}_1 \cdot a]_{\varphi_1}$, for $j = 1, \dots, \ell$; hence 
$\{{\overline q}_1 \cdot a, \ \ldots \ , {\overline q}_{\ell} \cdot a \}$
$\subseteq$ $[{\overline q}_1 \cdot a]_{\varphi_1}$.
Also, since ${\overline q}_j$ is ${\sf mod} \varphi_1$-equivalent to 
${\overline q}_1$
and since ${\sf mod} \varphi_1$ refines ${\sf mod} \varphi$, we have for
all $j$:
 \ $({\overline q}_j)\varphi = ({\overline q}_1)\varphi$ \ and \  
$({\overline q}_j \cdot a)\varphi = ({\overline q}_1 \cdot a)\varphi$.

In {\bf A} we consider the $\equiv_{\mathcal R}$-class 
$R = [({\overline q}_1 a)\varphi]_{\equiv_{\mathcal R}}$, and the 
corresponding $\equiv_{\mathcal R}$-class $\ov{R}$ $=$
$[\ov{q}_1 a]_{\equiv_{\mathcal R}}$ in ${\bf A}^{\sf cov}$.  
Since $({\overline q}_j)\varphi = ({\overline q}_1)\varphi$, we have
$R = [({\overline q}_1 a)\varphi]_{\equiv_{\mathcal R}} = $
$[({\overline q}_j a)\varphi]_{\equiv_{\mathcal R}}$ for
$j = 1, \ldots, \ell$.  
Recall that in ${\bf A}^{\sf cov}$ we define ${\sf rank}(\ov{R})$ to be
the maximum of the ranks of all the edges in $\ov{R}$, and that the rank of
an edge $(q_1, a, q_2)$ in $\ov{R}$ is defined to be the rank of the image 
edge $((q_1)\varphi, a, (q_2)\varphi)$ in $R$. 
Let $E(S)$ be the set of edges between states in a set $S$. Below, 
${\sf paths}(.)$ is taken in the $\cal R$-class $R$.

We now proceed in a number of stages. 

\medskip

\noindent {\sc Stage 1:} \ Let 

\smallskip

$n_1 \ = \ {\sf max}\{ \, {\sf rank}(\ov{e}) \ : \ \ov{e} \, \in \, $
$\bigcup_{1 \le j \le \ell}$
$E([{\overline q}_j a]_{\equiv_{\mathcal R}}) \ \subseteq \ $
$E({\bf A}^{\sf cov}) \, \}$ \ \ (where ${\sf max}(\varnothing) = -1$),

\medskip

$S_1 \ = \ {\sf paths}\big( (\ov{q}_1 a)\varphi, \ n_1 \big)$
  $ \ \cup \ \{\varepsilon\}$,

\medskip

$\ov{P}_1 \ = \ \{ \ov{q}_j a t \in Q^{\sf cov} \ : \ t \in S_1,$
  $ \ 1 \le j \le \ell \}$.

\medskip

\noindent {\sc Stage $i+1$:} \ Assuming $S_h, \, {\ov P}_h$ have been 
defined for $1 \leq h \leq i$, let  

\smallskip

$n_{i+1} \ = \ {\sf max}\{ \, {\sf rank}(\ov{e}) \ : \ \ov{e} \, $
is an edge of an $\cal R$-class of $\ov{P}_i \, \}$, 

\medskip

$S_{i+1} \ = \ {\sf paths}\big( (\ov{q}_1 a)\varphi, \ n_{i+1} \big)$
  $ \ \cup \ \ \{\varepsilon\}$,

\medskip

${\ov P}_{i+1} \ = \ \bigcup \, \{ [\ov{q}_j a t]_{\equiv_{\cal R}} \ : $
          $ \  \ov{q}_j a t \in Q^{\sf cov}, \ $
          $t \in S_{i+1}, \ 1 \le j \le \ell \}$,

\medskip

\noindent where $[\ov{q}_j a t]_{\equiv_{\cal R}}$ is the $\cal R$-class
of $\ov{q}_j a t$ in ${\bf A}^{\sf cov}$.

\medskip

\noindent {\sc Stage End:} \    
Continuing in this way we construct chains    

\smallskip

$n_1 \le n_2 \le \ \ldots \ \le n_i \le \ \ldots \ $
$ \le n_{\infty} = {\sf max} \{n_i : i = 1, 2,  \ \ldots \ \}$, 

\medskip

$S_1 \subseteq S_2 \subseteq \ \ldots \ \subseteq S_i \subseteq $
$ \ \ldots \ \subseteq \ S_{\infty} = \bigcup_i S_i \ \subseteq A^*$, 
 \ and \ 

\medskip

$\ov{P}_1 \subseteq \ov{P}_2 \subseteq \ \ldots \ \subseteq \ov{P}_i$ 
$\subseteq \ \ldots \ \subseteq \ $
$\ov{P}_{\infty} = \bigcup_i {\ov P}_i \ \subseteq Q^{\sf cov}$.

\medskip

\noindent These sequences are actually finite. Indeed, 
${\ov P}_i \subseteq Q^{\sf cov}$, which is a fixed finite set. Hence, 
$n_i$ is bounded since it is defined in terms of ${\ov P}_{i-1}$. Hence the 
sequence $S_i$ is of bounded, being defined in terms of $n_i$.

To define $\varphi_1$, $\varphi_2$, and ${\bf A}^{\sf cov \, GLC}$,
we want to construct the ${\sf mod} \varphi_1$ congruence class 
$[\ov{q}_j a]_{\varphi_1}$ so that 
$\ov{q}_j a \stackrel{\varphi_1}{\longmapsto} [\ov{q}_j a]_{\varphi_1}$
$\stackrel{\varphi_2}{\longmapsto} (\ov{q}_1 a)\varphi$, where 
$\varphi_1 \varphi_2 = \varphi$, and where $\ov{q}_j a$ and 
$(\ov{q}_j a)\varphi = (\ov{q}_1 a)\varphi$ are known. And we want 
$\varphi_2$ to be injective on each $\cal R$-class. So 
$([(\ov{q}_1 a)\varphi]_{\equiv_{\cal R}})\varphi^{-1}$ is a union of 
$\cal R$-classes of ${\bf A}^{\sf cov}$, each of which maps into  
$[(\ov{q}_1 a)\varphi]_{\equiv_{\cal R}}$.
By construction, $\ov{P}_{\infty}$ is a union of $\cal R$-classes of
${\bf A}^{\sf cov}$.

Therefore, for $\ov{q}_j a \in \ov{P}_{\infty}$ we define 
$[\ov{q}_j a]_{\varphi_1}$ to be  \     
$(\ov{q}_1 a)\varphi \varphi^{-1} \ \cap \ \ov{P}_{\infty}$.
More generally, for any $\ov{q} \in \ov{P}_{\infty}$ we define  
$[\ov{q}]_{\varphi_1}$ to be

\medskip

 \ \ \  \ \ \   
 $(\ov{q})\varphi \varphi^{-1} \ \cap \ \ov{P}_{\infty}$.

\medskip

\noindent Then $\varphi_2$ is defined by 

\medskip

 \ \ \  \ \ \ $[\ov{q}]_{\varphi_1}$ $ \ \longmapsto \ $
$([\ov{q}]_{\varphi_1})\varphi \ \in Q$  

\medskip

\noindent where $Q$ is the state set of {\bf A}. 
This completes the inductive step of the construction of 
${\bf A}^{\sf cov\, GLC}$.

\subsection{Proof of correctness}

\begin{pro} \  Suppose that {\bf A} is the Cayley graph of an 
$A$-monoid $M_A$ such that the stabilizers of ${\bf A}^{\sf cov}$ are
$\cal R$-trivial bands with unambiguous $\cal L$-order, and such that 
the $\cal R$-order of $M_A$ is a tree.
Then ${\sf glc} = {\sf GLC}$.  
\end{pro}
{\bf Proof.} \   
By construction $\varphi$, $\varphi_1$, $\varphi_2$ are $A$-automaton 
morphisms satisfying $\varphi = \varphi_1 \circ \varphi_2$.
We want to show that $\varphi_2: {\bf A}^{\sf cov \, GLC} \to {\bf A}$ is 
1:1${\mathcal R}$ and obeys the rank condition; and we want to show 
that $\varphi_2$ is maximal among all the right-factors of $\varphi$ that 
are 1:1${\mathcal R}$ and that satisfy the rank condition. 

\medskip

\noindent {\sf Proof that $\varphi_2$ is 1:1${\mathcal R}$ and
satisfies the rank condition}

\smallskip

Let ${\ov q} \in Q^{\sf cov}$. If $(\ov{q}a)\varphi$ is in the 
$\cal R$-class $R$, and $(\ov{q}a)\varphi_1$ is in the $\cal R$-class 
$\tilde{R}$, then by construction (see also Prop.\ \ref{edgesPaths}), 
$E\big((\ov{q}a)\varphi, \, n_{\infty}\big) \subseteq E(R)$;
recall that $E(r,j)$ denotes the set of edges of ${\sf paths}(r,j)$.
Then, clearly, $\varphi_2$ is 1:1$\cal R$ and satisfies the rank 
condition, since it embeds the edges of 
${\sf paths}\big((\ov{q}a)\varphi, \, n_{\infty}\big)$ into $E(R)$.

\medskip

\noindent {\sf Proof of maximality of $(.)^{\sf cov \, GLC}$ with
respect to the 1:1${\mathcal R}$ property and the rank condition }
 
\smallskip

We want to show that $\varphi_2$ is the maximal right-factor of $\varphi$ 
that is 1:1${\mathcal R}$ and satisfies the rank condition.  
Let {\bf B} be an $A$-automaton, and let 
$\theta_1: {\bf A}^{\sf cov} \twoheadrightarrow {\bf B}$ and 
$\theta_2: {\bf B} \twoheadrightarrow {\bf A}$ be automaton morphisms, where
$\theta_2$ is 1:1${\mathcal R}$ and satisfies the rank condition, and such
that $(.)\varphi = (.)\theta_1 \theta_2$. Hence the congruence 
${\sf mod} {\theta_1}$ refines ${\sf mod} {\varphi}$.  We want to show
that the congruence ${\sf mod} {\varphi_1}$ refines ${\sf mod} {\theta_1}$. 
  
\smallskip

\noindent {\sc Claim.} \ For any ${\ov q} \in A^{\sf cov}$,  
 \ $[\ov q]_{\varphi_1} \subseteq [\ov q]_{\theta_1}$.

\smallskip

\noindent Proof of the Claim: We follow the inductive construction of 
${\sf mod} {\varphi_1}$.  First, ${\ov q} \in [\ov q]_{\theta_1}$, 
obviously. 

We will use $E(V)$ to denote the set of edges between vertices in $V$ 
(including loops).
We have: \ $[{\ov q}]_{\equiv_{\mathcal R}} \subseteq $
$\bigcup_{{\ov p} \in [{\ov q}]_{\theta_1}} [{\ov p}]_{\equiv_{\mathcal R}}$
$=$ $[[{\ov q}]_{\theta_1}]_{\equiv_{\mathcal R}}$;
hence, since $\theta_2$ is 1:1$\mathcal R$, we obtain \ $E_1^o \subseteq $
$E\big(\big( [[{\ov q}]_{\theta_1}]_{\equiv_{\mathcal R}}\big)\varphi \big)$.
Since $\theta_2$ satisfies the rank condition we also have $E_1 \subseteq$
$E\big(\big( [[{\ov q}]_{\theta_1}]_{\equiv_{\mathcal R}}\big)\varphi \big)$.
Thus, ${\ov P}_1 \subseteq [{\ov q}]_{\theta_1}$.

Continuing along the inductive construction of ${\sf mod} {\varphi_1}$ we 
obtain ${\ov P}_{\infty} \subseteq [{\ov q}]_{\theta_1}$. Hence,
$[\ov q]_{\varphi_1} \subseteq [\ov q]_{\theta_1}$.
This proves the Claim.

\smallskip

Since the Claim holds for every state ${\ov p}$ in $[\ov q]_{\theta_1}$,
we conclude that $[\ov q]_{\theta_1}$ is a union of congruence classes
$[{\ov p}]_{\varphi_1}$.  \ \ \ $\Box$
 


\section{The Key Lemma}

We will use $\liteq$ to denote equality in $A^+$ or $A^*$, i.e., literal
equality of words. 

\begin{defn}
 \ A cover $(.)^{\sf cov}$ of $A$-automata has the {\em backwards-$k$}
property (for $k \geq 1$) \ iff \ for every $A$-automaton 
${\bf A} = (Q, \i, \cdot)$ and every $\alpha, \beta \in A^+$, the following
holds.  If $\alpha \beta = \alpha$ in the syntactic semigroup 
$S_{_A}^{\sf cov}$ of ${\bf A}_{_A}^{\sf cov}$, then $\alpha$ can be 
factored in $A^+$ as 
$\alpha \liteq \tilde{\alpha} \beta_1 \ldots \beta_k$ (with
$\tilde{\alpha}, \beta_1, \ \ldots, \beta_k \in A^+$) such that in 
the syntactic semigroup $S_A$ of {\bf A} we have $\tilde{\alpha} = \alpha$ 
and $\beta_1 = \beta_2 = \ \ldots \ = \beta_k = \beta = \beta^2$.
\end{defn}

\begin{pro} \label{backkExp} 
 \ For $A$-generated semigroups, let $(.)_A^{\sf IS}$ be an expansion with
$\cal R$-trivial idempotent stabilizers, let $(.)_A^{\wedge_{(k+1)}}$ 
be the Henckell $(k+1)$-factor expansion, and let $(.)_A^{\sf RB}$
be the rectangular-bands expansion.
Then for any $A$-generated semigroup $S$, the triply expanded semigroups 
$S_{\, _A}^{{\sf IS} \, \wedge_{(k+1)} \, {\sf RB}}$ and
$S_{\, _A}^{ \wedge_{(k+1)} \, {\sf IS} \ {\sf RB}}$ 
have the backwards-$k$ property.  
\end{pro}
{\bf Proof.} If $\alpha \beta = \alpha$ holds in 
$S_{\, _A}^{{\sf IS} \, \wedge_{(k+1)} \, {\sf RB}}$ or in
$S_{\, _A}^{ \wedge_{(k+1)} \, {\sf IS} \ {\sf RB}}$ then 
$\alpha \beta = \alpha$ also holds in 
$S_{\, _A}^{{\sf IS} \, \wedge_{(k+1)}}$, respectively
$S_{\, _A}^{ \wedge_{(k+1)}}$. 
The existence of a factorization 
$\alpha = \tilde{\alpha} \beta_1 \ldots \beta_k$ in $S_A$ then follows from
the basic properties of the expansion $(.)_A^{\wedge_{(k+1)}}$.  
 \ \ \ $\Box$

\bigskip

For an $A$-automaton {\bf A} and an expansion $(.)^{\sf exp}$, let 
${\bf A}^{\sf exp \, GLC}$ be the {\sf glc}-cover in the interval 
$[{\bf A}, \, {\bf A}^{\sf exp}]$, with start state $\tilde{\i}$ and state 
set $Q_{\sf glc}$. From ${\bf A}^{\sf exp \, GLC}$ and a string 
$s = a_1 a_2 \ \ldots \ a_n$, where $a_i \in A$, we define 
a {\em string automaton} ${\sf str}(s)$ as follows.
In the state graph of ${\bf A}^{\sf exp \, GLC}$ we consider any walk 

\smallskip

\hspace{.3in}
$\tilde{\i} \stackrel{a_1}{\longrightarrow} u_1$
$\stackrel{*}{\longrightarrow} v_1$
$\stackrel{a_2}{\longrightarrow} \ \ \ldots \ \ $
$\stackrel{a_i}{\longrightarrow} u_i$
$ \stackrel{*}{\longrightarrow} v_i$
$\stackrel{a_{i+1}}{\longrightarrow}  \ \ \ldots \ \ $
$\stackrel{a_n}{\longrightarrow} u_n \in R_n$,

\smallskip

\noindent where 
$u_1, v_1, \ \ldots, \ u_{n-1}, v_{n-1}, u_n \in Q_{\sf glc}$ are such 
that $u_i \equiv_{\cal R} v_i$ for $i = 1, \ \ldots, \ n-1$; there is no 
$v_n$.
Let $R_i$ be the reachability class of $u_i$ (and of $v_i$) in 
${\bf A}^{\sf exp \, GLC}$. 
We assume $R_i \neq R_j$ (hence they are vertex-disjoint) when $i \neq j$; 
hence $R_i >_{\cal R} R_{i+1}$ for all $i$.
The state set of ${\sf str}(s)$ is
 \ $\{ \tilde{\i} \} \ \cup \ \bigcup_{1 \le i \le n} R_i$,
and the arrows of ${\sf str}(s)$ are just the 
${\bf A}^{\sf exp \, GLC}$-arrows between states of ${\sf str}(s)$.
We will use $\tilde{\i}$ as the start state, and $R_n$ as the set of accept
states of ${\sf str}(s)$.
Conversely, ${\sf str}(s)$ uniquely determines the above walk: $u_i$ is 
the entry point into $R_i$, and $v_i$ is the exit point.  

Any $A$-generated monoid $M_A$ will be viewed as an $A$-automaton via
its right Cayley graph.

\begin{thm} {\bf (``Key Lemma'').} \label{KL}
Let $M_A$ be an $A$-generated monoid, viewed as an $A$-automaton {\bf A}.
We consider the $A$-generated expansion 
$M_{\, _A}^{\wedge_{(k+1)} \, {\sf IS} \ {\sf RB}}$ (for some
$k >1$), viewed as an $A$-automaton ${\bf A}^{\sf exp}$.
Let ${\bf A}^{\sf exp \, GLC}$ be the {\sf glc}-cover of $A$-automata in
the interval $[{\bf A}, \, {\bf A}^{\sf exp}]$.
Let ${\sf str}(a_1 a_2 \ldots a_n)$ be a string automaton in
${\bf A}^{\sf exp \, GLC}$, with $\cal R$-classes $R_i$ for 
$i = 1, \ \ldots, \ n$. Finally, let 
$r_n = {\sf max}\{ {\sf rank}_{\bf A}(e) : e \in E(R_n)\}$, where 
$E(R_n)$ is the edge-set of the reachability class $R_n$.

Then there exist $d_1, \ \ldots, \ d_{n-1}, d_n \in A^*$ such that $d_i$ 
labels a path $u_i \stackrel{d_i}{\longrightarrow} v_i$ in $R_i$ (in 
${\bf A}^{\sf exp \, GLC}$), $1 \leq i \leq n-1$, and such that the word 
$w \liteq a_1 d_1 a_2 \ \ldots \ a_{n-1} d_{n-1} a_n d_n$ 
satisfies:

\smallskip

\noindent {\bf (1)} \ \ There is $\ell_n \in A^*$ such that in 
${\bf A}^{\sf exp}$: \ $w \, \ell_n = w$.

\smallskip

\noindent {\bf (2)}  \ \    
${\sf rank}_{\bf A}(\ell_n) \ = \ r_n \ > \ {\sf rank}_{\bf A}(d_n)$.

\smallskip

\noindent {\bf (3)} \  {\rm (Backwards-$k$ property):} \ There exist 
$w', t_1, \ldots, \ t_k \in A^+$ (where $k$ is as in $\wedge_{(k+1)}$) 
such that $w \liteq w' \, t_k \ \ldots \, t_2 \, t_1$ in $A^+$, 
and in $M_A$: \ $w = w'$, \ $t_k = \ \ldots \ = t_1 = \ell_n$.

The path labeled by $t_1$, when read backwards from the state 
${\tilde q}_n$ that $\ell_n$ loops on, eventually visits the source 
state $v_{n-1}$ of the edge labeled by $a_n$ in 
${\sf str}(a_1 a_2 \ldots a_n)$; i.e., $a_n d_n$ is a suffix of $t_1$,
as seen in the next graph $ \ \ \ \ldots \ $ 
$\stackrel{*}{\longrightarrow} v_{n-1}$
$\stackrel{a_n}{\longrightarrow} u_n$
$\stackrel{d_n}{\longrightarrow} {\tilde q}_n$
$\circlearrowleft \ell_n$ \, (in ${\bf A}^{\sf exp \, GLC}$).  
\end{thm}
{\bf Proof.} Let us pick $d_1, \ \ldots, \ d_{n-1}$ arbitrarily so that 
$u_i \stackrel{d_i}{\longrightarrow} v_i$ in $R_i$, for $1 \leq i \leq n-1$. 
Reading $w_1 \liteq a_1 d_1 a_2 \ \ldots \ a_{n-1} d_{n-1} a_n$ from the 
start state $\ov{\i}$ in ${\bf A}^{\sf exp}$, we reach the state 
$\ov{\i} \cdot w_1$, which maps into $R_n$ by the map 
${\bf A}^{\sf exp} \twoheadrightarrow {\bf A}^{\sf exp \, GLC}$.
But $\ov{\i} \cdot w_2$ does not map into $R_n$, where $w_2 = w_1 a_n^{-1}$.
(We use the notation $wa a^{-1} = w$; i.e., $a^{-1}$ is the operation of
removing the last letter $a$ from a word that ends in $a$, and $a^{-1}$ is
undefined on words that do not end in $a$.)

Thus we can run {\sc Stage} 1 of the Bottom-up construction on state
${\ov q} \cdot a$, where for ${\ov q}$ we take ${\ov \i} \cdot w_2$ in 
${\bf A}^{\sf exp}$, and for $a$ we take $a_n$. 
Then we run the inductive {\sc Stage} $i+1$ for increasing $i$ until we 
reach $n_{\infty}$. 

If $n_1 = n_{\infty}$ then as we enter the $\cal R$-class 
$[{\ov q} a]_{\equiv_{\cal R}}$ we have maximum edge-rank equal to 
$n_{\infty}$; we let $d_n$ be the empty word, and let $\ell_n$ be a loop
at the entrance of $[{\ov q} a]_{\equiv_{\cal R}}$, passing through all the
edges in $E([{\ov q} a]_{\equiv_{\cal R}})$.

If $n_1 < n_{\infty}$ we go to the 2nd stage of the induction, i.e.,  
we find an $\cal R$-class $\ov R$ of ${\bf A}^{\sf exp}$ with rank
$n_2 > n_1$.  So we have a path 
 \ $\ldots \ \stackrel{d_n}{\longrightarrow} {\ov R}_2$,
and this path up to (but excluding) ${\ov R}_2$ has edges of rank $n_1$. 
Then we take $d_n$ as indicated, and we take $\ell_n$ to be a loop at the 
entrance of ${\ov R}_2$, passing through all the edges in $E({\ov R}_2)$.

We continue with $n_3$, etc., until $n_{\infty}$ is first encountered.
Now, item (3) follows from Prop.\ \ref{backkExp}, applied to the 
$\ell_n$-loop in \ \ $\ldots$
$ \stackrel{d_n}{\longrightarrow} \bullet \circlearrowleft \ell_n$.
 \ \ \   \ \ \ $\Box$

\bigskip

The Key Lemma can be generalized by replacing $k$ by two parameters,
$k_1$ and $k_2$. For this we use the expansion \   
$(.)_{\, _A}$\hspace{-.1in}
$^{\wedge_{k_1} \, {\sf RB} \ \wedge_{(k_2 + 1)} \ {\sf IS} \ {\sf RB}}$.

\begin{thm} {\bf (Key Lemma with $k_1$ and $k_2$).} \label{KLk1k2} \  
Under the assumptions of Theorem \ref{KL}, except that the expansion is
now based on \
$(.)_{\, _A}$\hspace{-.1in}
$^{\wedge_{k_1} \, {\sf RB} \ \wedge_{(k_2 + 1)} \ {\sf IS} \ {\sf RB}}$, 
we have the same conclusions, except for changes in item (3):

\smallskip

\noindent {\bf (3)} \ {\rm (Backwards-$k_1$-$k_2$ property):} There exist
$w', t_1, \ldots, \ t_{k_2} \in A^+$ such that 
$w \liteq w' \, t_{k_2} \ \ldots \, t_2 \, t_1$ in $A^+$,
and in $M_A^{\wedge_{k_1} \, {\sf RB}}$: \   
$w = w'$, \ $t_{k_2} = \ \ldots \ = t_1 = \ell_n$.

The path labeled by $t_1$, when read backwards from the state 
${\tilde q}_n$ that $\ell_n$ loops on, eventually visits the source
state $v_{n-1}$ of the edge labeled by $a_n$ in
${\sf str}(a_1 a_2 \ldots a_n)$; i.e., $a_n d_n$ is a suffix of $t_1$, 
as seen in the next graph $ \ \ \ \ldots \ $
$\stackrel{*}{\longrightarrow} v_{n-1}$
$\stackrel{a_n}{\longrightarrow} u_n$
$\stackrel{d_n}{\longrightarrow} {\tilde q}_n$
$\circlearrowleft \ell_n$ \, (in ${\bf A}^{\sf exp \, GLC}$).

\end{thm}
{\bf Proof.} The proof is the same as for Theorem \ref{KL}, but we use
the morphisms \  
$M_A^{\wedge_{k_1} \, {\sf RB} \ \wedge_{(k_2 + 1)} \ {\sf IS} \ {\sf RB}}$
$\twoheadrightarrow$ $M_A^{\wedge_{k_1} \, {\sf RB}}$
$\twoheadrightarrow$ $M_A^{\wedge_{k_1}}$  for item (3). 
 \ \ \ $\Box$

\bigskip

\noindent {\bf Notation:}  Let $({\sf GM}, A)$ be the starting {\em group
mapping} semigroup whose complexity we are trying to compute, let 
$M_A = ({\sf GM}, A)^{\sf IS \, RB}$, and let $(.)^{\sf exp}$ be
$(.)^{\wedge_{k_1}\, {\sf RB}\ \wedge_{(k_2 + 1)}\ {\sf IS}\ {\sf RB}}$.  
Then the semigroup automaton ${\bf A}_A^{\sf exp \ GLC}$ will be called 
{\sc PreFF}$(k_1, k_2)$ (for ``pre-funny fractal''). 

\smallskip

For any finite semigroup $S$, let $\omega(S)$ be the smallest positive 
integer $m$ such that for all $s \in S$, $s^m$ is an idempotent.  

%
%
%
%
%
%

\bigskip 

{\it To be continued.}


\bigskip

\bigskip


\end{document}